\documentclass[11pt]{amsart}
\usepackage{amssymb,amsmath, amsthm, amsfonts}

\usepackage[usenames,dvipsnames]{xcolor}
\usepackage{graphicx}
\usepackage{listings}

\usepackage[normalem]{ulem}

\usepackage[margin=1in]{geometry}
\usepackage{lstautogobble}
\usepackage{enumerate}
\usepackage[shortlabels]{enumitem}
\usepackage{thmtools}
\usepackage{thm-restate}
\usepackage{amsthm}
\usepackage{verbatim}
\usepackage{accents}
\usepackage{mathtools}
\usepackage{physics}

\usepackage[colorlinks=true, allcolors=magenta,linkcolor=blue, citecolor=magenta]{hyperref}
\usepackage[capitalise]{cleveref}
\usepackage[bottom]{footmisc}
\crefformat{equation}{(#2#1#3)}

\usepackage{float}
\restylefloat{table}

\usepackage{mathrsfs}
\setlist{  
	listparindent=\parindent,
	parsep=0pt,
}

\theoremstyle{plain}
\newtheorem{thm}{Theorem}[section]
\newtheorem{prop}[thm]{Proposition}

\newtheorem{cor}[thm]{Corollary}

\theoremstyle{definition}

\newtheorem{remark}[thm]{Remark}

\Crefname{thm}{Theorem}{Theorems}
\Crefname{prop}{Proposition}{Propositions}

\numberwithin{equation}{section} 


\DeclarePairedDelimiter\ipp{\langle}{\rangle}

\DeclarePairedDelimiter{\pa}{\lparen}{\rparen}

\DeclarePairedDelimiter{\jp}{\langle}{\rangle}

\DeclareMathOperator{\supp}{supp}

\DeclareMathOperator{\dist}{dist}
\DeclareMathOperator{\sgn}{sgn}

\newcommand{\M}{{\mathcal{M}}}

\newcommand{\p}{{\partial}}

\renewcommand{\d}{\mathsf{d}}

\newcommand{\R}{{\mathbb{R}}}

\newcommand{\N}{{\mathbb{N}}}
\newcommand{\Q}{{\mathbb{Q}}}

\newcommand{\Ss}{{\mathbb{S}}}

\newcommand{\T}{{\mathbb{T}}}
\newcommand{\g}{{\mathsf{g}}}
\newcommand{\G}{{\mathsf{g}}}

\newcommand{\Sc}{{\mathcal{S}}}

\renewcommand{\M}{{\mathbb{M}}}
\newcommand{\I}{\mathbb{I}}

\renewcommand{\k}{\mathsf{k}}
\newcommand{\ga}{\gamma}
\newcommand{\nab}{\nabla}

\newcommand{\tl}{\tilde}

\newcommand{\nn}{\nonumber}

\newcommand{\XN}{X_N}

\newcommand{\ux}{X}

\newcommand{\uz}{Z}
\newcommand{\ep}{\epsilon}
\newcommand{\vep}{\varepsilon}
\newcommand{\al}{\alpha}
\newcommand{\be}{\beta}
\newcommand{\ka}{\kappa}
\newcommand{\la}{\lambda}
\newcommand{\Om}{\Omega}

\newcommand{\indic}{\mathbf{1}}

\newcommand{\Fr}{\mathsf{F}}
\newcommand{\Hr}{\mathsf{H}}

\newcommand{\Uu}{\mathfrak{U}}
\newcommand{\As}{\mathsf{A}}

\newcommand{\E}{{\mathbb{E}}}

\newcommand{\Dm}{|\nabla|}

\newcommand{\rs}{\mathsf{r}}
\newcommand{\cd}{\mathsf{c}_{\mathsf{d},\mathsf{s}}}

\newcommand{\s}{\mathsf{s}}

\renewcommand{\k}{\mathsf{k}}

\newcommand{\zg}{|z|^\ga}

\renewcommand{\P}{\mathcal{P}}

\newcommand{\as}{\mathsf{a}}
\let\div\relax
\DeclareMathOperator{\div}{\mathrm{div}}

\def\XXint#1#2#3{{\setbox0=\hbox{$#1{#2#3}{\int}$ }
		\vcenter{\hbox{$#2#3$ }}\kern-.6\wd0}}

\let\oldtocsection=\tocsection

\let\oldtocsubsection=\tocsubsection

\let\oldtocsubsubsection=\tocsubsubsection

\renewcommand{\tocsection}[2]{\hspace{0em}\oldtocsection{#1}{#2}}
\renewcommand{\tocsubsection}[2]{\hspace{1em}\oldtocsubsection{#1}{#2}}
\renewcommand{\tocsubsubsection}[2]{\hspace{2em}\oldtocsubsubsection{#1}{#2}}

\title[Commutators, mean-field, and supercritical mean-field limits]{Commutators, mean-field, and supercritical mean-field limits for Coulomb/Riesz gases}

\author[M. Rosenzweig]{Matthew Rosenzweig}
\address{Matthew Rosenzweig, Carnegie Mellon University, Department of Mathematical Sciences, Pittsburgh, PA} 
\email{mrosenz2@andrew.cmu.edu}
\thanks{M.R. was supported by NSF grants DMS-2441170, DMS-2345533, DMS-2342349.}

\begin{document}
 
	\begin{abstract} 
	This paper is intended as a companion to the author's talk ``Commutator estimates and mean-field limits for Coulomb/Riesz gases'' at the 2025 \emph{Journ\'ees \'equations aux d\'eriv\'ees partielles} in Aussois. The goal is to provide a concise, accessible account of  sharp commutator estimates recently obtained for modulated energies associated to Coulomb/Riesz interactions and how these estimates lead to optimal results for mean-field and supercritical mean-field limits of Coulomb/Riesz gas dynamics via the modulated-energy method. The exposition centers on the works \cite{rosenzweig_sharp_nodate,rosenzweig_lake_2025} with Serfaty and \cite{hess-childs_sharp_2025,hess-childs_sharp_2025} with Hess-Childs and Serfaty.  
	\end{abstract}
	\maketitle

\section{Introduction}\label{hofis:sec:intro}
When studying systems of $N$ distinct points $\ux_N = (x_1, \dots, x_N)\in (\R^\d)^N$  with interaction energy  
\begin{equation}
\sum_{1\le i\neq j\le N} \g(x_i,x_j),
\end{equation}
one is led to comparing the sequence of empirical measures $\mu_N \coloneqq \frac1N \sum_{i=1}^N \delta_{x_i}$ to an  average, or \emph{mean-field}, density $\mu$. This comparison is conveniently performed by considering a {\it modulated energy}, or (squared) ``distance''  between $\mu_N$ and $\mu$, defined by 
\begin{equation}\label{hofis:eq:modenergy}
\Fr_N(\ux_N, \mu) \coloneqq \frac12\int_{(\R^\d)^2 \setminus \triangle} \g(x,y) d\Big(\frac{1}{N} \sum_{i=1}^N \delta_{x_i} - \mu\Big)(x)d\Big(\frac{1}{N} \sum_{i=1}^N \delta_{x_i} - \mu\Big)(y),
\end{equation}
where $\triangle\subset(\R^\d)^2$ denotes the diagonal. Physically, the quantity corresponds to the total energy of a system of $N$ discrete charges located at the positions $x_i$ neutralized by a background charge distribution $\mu$ with the infinite self-interaction of each charge removed, the latter corresponding to the excision of the diagonal.

{The quantity \eqref{hofis:eq:modenergy} originated in the study of the statistical mechanics of Coulomb/Riesz gases \cite{sandier_1d_2015,sandier_2d_2015,rougerie_higher-dimensional_2016, petrache_next_2017} and was later used in the derivation of mean-field dynamics \cite{duerinckx_mean-field_2016,serfaty_mean_2020,nguyen_mean-field_2022} and following works. The notion of modulated energy to prove a stability estimate is classical, with the term first appearing to this author's knowledge in the influential work of Brenier \cite{Brenier2000}. In the nonsingular case, where the diagonal may be reinserted, the modulated energy coincides with the square of what is known in the statistics literature as \emph{maximum mean discrepancy (MMD)} \cite{gretton_kernel_2006, gretton_kernel_2007,gretton_kernel_2012}, a type of integral probability metric \cite{muller_integral_1997}. MMDs for varying choices of kernels are widely used in statistics and machine learning contexts as distances between probability measures (see, e.g., \cite{modeste_characterization_2024, kolouri_generalized_2022,  hertrich_generative_2023}).} 

An essential point in this comparison is to control derivatives of $\Fr_N$ along  a transport $v:\R^\d\rightarrow\R^\d$, i.e. the quantities
\begin{equation}\label{15}
 \frac{d^n}{dt^n}\Big|_{t=0} \Fr_N( (\I + tv)^{\oplus N} (\ux_N), (\I + tv)\# \mu)= {\frac12}\int_{(\R^\d)^2\setminus \triangle} 
\nabla_x^{\otimes n} \g(x,y):  (v(x)-v(y))^{\otimes n}  d ( \mu_N- \mu)^{\otimes 2},
\end{equation}
where $\I:\R^\d\rightarrow\R^\d$ is the identity, $(\I+t v)^{\oplus N} (\ux_N) \coloneqq (x_1 + tv(x_1), \ldots, x_N+ tv(x_N))$, and $:$ denotes the inner product between the tensors. Here, $\nab^{\otimes n}$ denotes the $n$-tensor with components $(\p_{\vec{\al}})_{\vec{\al}\in [\d]^n}$. 

For applications to the mean-field limit, $v$ is the velocity field of the limiting dynamics as $N \to \infty$. For applications to central limit theorems (CLTs) for the fluctuations (the next-order description), $v$ is the gradient of a test function evolved along the adjoint linearized mean-field flow \cite{huang_fluctuations_nodate,cecchinConvergenceRateFluctuations2025}, and these inequalities are at the core of the ``transport'' approach to fluctuations of canonical Gibbs ensembles \cite{leble_fluctuations_2018,BLS2018,serfaty_gaussian_2023,peilen_local_2025} (cf. loop/Dyson--Schwinger equations in random matrix theory, e.g.  \cite{BG2013, BG2024, BBNY2019}).

The desired control is a functional inequality of the form
\begin{equation}\label{eq:introcomm}
|\eqref{15} | \le C_1(\Fr_N(\ux_N, \mu) + C_2N^{-\alpha}),
\end{equation}
where $\alpha \in (0,\infty]$ depends on $\d,\g$, $C_1>0$ is a constant depending on $\d,\g,n,v$, and $C_2>0$ depends on $\d,n,\g,\mu$. As observed in \cite{rosenzweig_mean-field_2020},  the expressions \eqref{15}  may be viewed as the quadratic form of a \textit{commutator}, akin to the famous Calder\'{o}n commutator \cite{calderon_commutators_1980, coifman_commutateurs_1978, christ_polynomial_1987, seeger_multilinear_2019}. Indeed, ignoring the exclusion of the diagonal,
\begin{equation}
    \int_{(\R^\d)^2} (v(x)-v(y))\cdot\nabla_x\g(x,y)f(x)g(y)dxdy=\Big\langle f ,\Big[v\cdot,\nabla \g\ast\Big] g \Big\rangle_{L^2},
\end{equation}
where $\ipp{\cdot,\cdot}_{L^2}$ denotes the $L^2$ inner product. For this reason, estimates of the form \eqref{eq:introcomm} are called \emph{commutator estimates}. The expressions \eqref{15} also have a \emph{stress-energy tensor} structure, a point on which we elaborate in \cref{ssec:rcomsupC}. 

The discussion so far has concerned a general (conditionally) positive definite kernel $\g$. To have any hope of showing a functional inequality of the form \eqref{eq:introcomm}, one needs to impose some assumptions on $\g$. Of particular interest and relevance for applications are singular interactions, model examples of which are the \emph{log/Riesz} potentials\footnote{The log case being the $\s\rightarrow 0$ limit of the Riesz case, the term ``Riesz'' will always include the log case.} given by   
	\begin{equation}\label{eq:gmod}
		\g(x,y) = \g(x-y)\coloneqq \begin{cases}  \frac{1}{\s} |x-y|^{-\s}, \quad  & \s \neq 0\\
			-\log |x-y|, \quad & \s=0.
		\end{cases}
	\end{equation}
We  consider any dimension $\d\ge 1$ and assume that $-2<\s<\d$. Up to a normalizing constant $\cd$, the potential $\g$ is the fundamental solution of the fractional Laplacian $(-\Delta)^{\frac{\d-\s}{2}} = |\nab|^{\d-\s}$, i.e. $|\nabla|^{\d-\s}\g = \cd \delta_0$. The cases $0\le \s<\d$ and $-2<\s<0$ are respectively referred to as the \emph{singular} and \emph{nonsingular} regimes because of the behavior of $\g$ at the origin. The restriction $\s>-2$ is natural because the Riesz potential ceases to be conditionally positive definite at $\s=-2$.\footnote{{We say that a kernel $k:\R^\d\times\R^\d\rightarrow \R$ is \emph{conditionally positive definite (CPD)} iff for any signed Borel measure $\rho$ with zero mass, it holds that $\int_{(\R^\d)^2}k(x,y)d\rho(x)d\rho(y)\ge 0$.}\label{fn:cpd_def}} The restriction to the \emph{potential} regime $\s<\d$, in which $\g$ is locally integrable, is to exclude the \emph{hypersingular} case, which is not of the mean-field type considered in this paper (e.g. see \cite{hardin_large_2018,hardin_dynamics_2021}). The case $\s=\d-2$ corresponds to the classical \emph{Coulomb} potential from electrostatics/gravitation. Based on this threshold, it is convenient to call the interaction \emph{sub-Coulomb} if $\s<\d-2$ and \emph{super-Coulomb} if $\s>\d-2$. More generally, long-range interactions of the form \eqref{eq:gmod} play an important role in physics \cite{dauxois_dynamics_2002}, approximation theory \cite{borodachov_discrete_2019}, and machine learning \cite{altekruger_neural_2023,hertrich_generative_2023, hertrich_wasserstein_2023, hagemann_posterior_2023}---just to list a few areas. Motivations for considering systems of the form \eqref{eq:MFode} are extensively reviewed in the lecture notes \cite{serfaty_lectures_2024}. 

Specializing to Riesz potentials and focusing on the first-order $n=1$ case, the control \eqref{eq:introcomm} was first proved  by Lebl\'{e} and Serfaty \cite{leble_fluctuations_2018} in the 2D Coulomb case, then generalized by Serfaty to the super-Coulomb case  \cite{serfaty_mean_2020}. The proof relied on the identification of the aforementioned {stress-energy tensor} structure in \eqref{15} and using integration by parts. Subsequently, the author  \cite{rosenzweig_mean-field_2020} used the commutator structure in  \eqref{15} to give a new proof in the Coulomb case. This commutator perspective was then used together with Q.H. Nguyen and Serfaty \cite{nguyen_mean-field_2022} to show the desired functional inequality for the full singular regime as well as a broader class of Riesz-type interactions including Lennard-Jones-type potentials. The argument of \cite{nguyen_mean-field_2022} has also been extended to the nonsingular regime in the joint work \cite{rosenzweig_wasserstein_nodate} with Slep\^{c}ev and Wang. As recognized in \cite{rosenzweig_sharp_nodate}, the stress-energy and commutator perspectives are in fact two sides of the same coin. These functional inequalities are crucial not only for proving CLTs for fluctuations of Riesz gases \cite{leble_fluctuations_2018, serfaty_gaussian_2023, peilen_local_2025}, but also for deriving mean-field limits \cite{serfaty_mean_2020,rosenzweig_mean-field_2022-1,rosenzweigMeanfieldApproximationHigherdimensional2022,nguyen_mean-field_2022,chodron_de_courcel_sharp_2023,rosenzweig_sharp_nodate, porat_singular_2025} and large deviation principles \cite{hess-childs_large_2023}. They have also been used to show the joint classical and mean-field limits for quantum systems of particles \cite{GP2022},  and the supercritical mean-field limits of classical \cite{han-kwan_newtons_2021, rosenzweig_rigorous_2023, menard_mean-field_2024, rosenzweig_lake_2025} and quantum systems of particles \cite{rosenzweig_quantum_2021, porat_derivation_2023}. The reader will see two applications, mean-field and supercritical mean-field limits, below.

The main problem for these functional inequalities has been to determine the optimal size of the additive error, i.e.~the exponent $\alpha$ in \eqref{eq:introcomm}, which a priori depends on $\d$ and $\s$. In formulating optimality, the appropriate comparison is with the minimal size of $\Fr_N(\ux_N,\mu)$. When $0\le \s<\d$, it is known that $\Fr(\ux_N,\mu) + \frac{\log(N\|\mu\|_{L^\infty})}{2\d N}\indic_{\s=0}\ge -C\|\mu\|_{L^\infty}^{\frac\s\d}N^{\frac{\s}{\d}-1}$ for $C=C(\d,\s)>0$ and this lower bound is sharp \cite{rosenzweig_sharp_nodate, hess-childs_sharp_2025}. {The log term in the $\s=0$ case is \emph{not} an additive error like $\|\mu\|_{L^\infty}^{\frac\s\d}N^{\frac{\s}{\d}-1}$, but rather the right quantity to consider is $\Fr(\ux_N,\mu) + \frac{\log(N\|\mu\|_{L^\infty})}{2\d N}\indic_{\s=0}$.} In controlling the magnitude of \eqref{15}, one needs a right-hand side which is also nonnegative. Thus, the best error one may hope for is of size $N^{\frac{\s}{\d}-1}$. When $-2<\s<0$, the modulated energy is nonnegative (it is a squared MMD as mentioned above), and its minimal value is $\propto N^{\frac\s\d-1}$ \cite{hess-childs_optimal_nodate}. Moreover, the functional inequality holds without additive error because the interaction is nonsingular \cite{rosenzweig_wasserstein_nodate}. 

{Only very recently has the program to show the sharp $N^{\frac{\s}{\d}-1}$ error rate been completed by the author and collaborators: first, the Coulomb case \cite{leble_fluctuations_2018,serfaty_gaussian_2023, rosenzweig_rigorous_2023}, then the entire (super-)Coulomb case \cite{rosenzweig_sharp_nodate}, and finally the entire singular Riesz case \cite{hess-childs_sharp_2025}. The completion of this program has yielded the sharp rates of convergence for the mean-field limit in the modulated energy distance, as well as rates for supercritical mean-field limits under optimal scaling assumptions. It has also been crucial for studying fluctuations around the mean-field limit in and out of thermal equilibrium \cite{leble_fluctuations_2018,serfaty_gaussian_2023,peilen_local_2025,huang_fluctuations_nodate, rosenzweig_commutator_nodate, rosenzweig_cumulants_nodate}.}

In the remainder of this paper, we describe the results and applications advertised in the preceding discussion. The body is organized as follows. {In \cref{sec:rcom}, we present our sharp {commutator} estimates in the form of \cref{thm:FIprime}, discuss their proof (\Cref{ssec:rcomsupC,ssec:rcomsubC}), followed by localized and higher-order variants \Cref{hofis:thm:mainunloc,hofis:thm:FI} (\cref{ssec:rcomloc}), and conclude with a discussion of transport regularity and defective estimates (\cref{ssec:rcomtransreg}). In \cref{sec:MF}, we turn to mean-field limits, recalling the background and main result \cref{thm:mainMF} (\Cref{ssec:MFbackground,ssec:MFmain}), and then explaining the role of commutator estimates and modulated energy in the proof of \cref{thm:mainMF} (\cref{ssec:MFcomm}), including the regularity assumptions (\cref{ssec:MFreg}), and end with some comments on the positive-temperature generalization known as modulated free energy (\cref{ssec:MFposT}). Finally, \cref{sec:Lake} is devoted to supercritical mean-field limits: we review the combined mean-field/quasineutral regime (\cref{ssec:Lakeback}); discuss scaling assumptions and prior work (\cref{lake:ssec:introPW}); present an informal version of the main result \cref{lake:thm:mainSMF} (\cref{lake:ssec:introMPf}); outline the method of proof; and close with comments on {optimality} (\cref{lake:ssec:1DCou}) and extensions to the sub-Coulomb (\cref{lake:ssec:subCou}) and non-monokinetic regimes (\cref{lake:ssec:nonmono}).}

\bigskip

{\noindent\textbf{Acknowledgments.} The author thanks the organizers of the 2025 {\it Journ\'ees \'equations aux d\'eriv\'ees partielles} for the opportunity to present his research in the stimulating and scenic environment of Aussois, and he thanks the staff of the Centre Paul-Langevin for their gracious hospitality. He also thanks the Institute for Computational and Experimental Research in Mathematics (ICERM) for its hospitality, where part of the research for this project was carried out during the Fall 2021 semester program ``Hamiltonian Methods in Dispersive and Wave Evolution Equations,'' as well as the Courant Institute of Mathematical Sciences at NYU for their hospitality during his visits in November 2021 and April 2024, where {part} of the reported research was carried out. Finally, and most importantly, the author {thanks his} collaborators Sylvia Serfaty and Elias Hess-Childs, without whom the works described here would not have {come to} fruition.} 

\section{Commutator estimates}\label{sec:rcom}
We elaborate on the commutator estimates for Riesz potentials highlighted in the introduction. To compactify the notation, given a point configuration $\XN$, reference measure $\mu$, and vector field $v$, let us write
\begin{align}\label{eq:A1def}
    n\in\N, \qquad    \As_n[\XN,\mu,v] \coloneqq \int_{(\R^\d)^2\setminus\triangle}\nab\g(x-y)\cdot(v(x)-v(y))d\Big(\frac1N\sum_{i=1}^N\delta_{x_i}-\mu\Big)^{\otimes 2}(x,y).
\end{align}

Summarizing the state of the art for global, first-order commutator estimates is the following theorem, combining the results of \cite{nguyen_mean-field_2022,rosenzweig_global--time_2023, rosenzweig_sharp_nodate,hess-childs_sharp_2025, rosenzweig_wasserstein_nodate}. In particular, taking $p=\infty$ in \eqref{eq:FIsupC}, \eqref{eq:FIsubC1} below yields the announced sharp $N^{\frac\s\d-1}$ additive error. As commented above, there is no additive error in the nonsingular case \eqref{eq:FInonsing}. Comparing the sub-Coulomb estimates \eqref{eq:FIsubC1} and \eqref{eq:FIsubC2}, the transport regularity in the latter is optimal even though its additive error is not, while the additive error in the former is optimal {even though} the transport regularity is not. 
We discuss this point more in \cref{ssec:rcomtransreg}. Later, in \cref{ssec:rcomloc}, we discuss generalizations of \cref{thm:FIprime} to localized and higher-order commutators.

    \begin{thm}\label{thm:FIprime}
    Let $-2<\s<\d$. Let $\mu \in L^1\cap L^p$ for $\frac{\d}{\d-\s}<p\le \infty$ with $\int_{\R^\d}d\mu = 1$, and associated to $\mu$, define the length scales
    \begin{align}
        \la &\coloneqq (N\|\mu\|_{L^p})^{-\frac{1}{\d}}\label{eq:ladef},\\\
        \ka &\coloneqq (N^{\frac{1}{\s+1}}\|\mu\|_{L^p})^{-\frac{1}{\d}}, \label{eq:kadef}
    \end{align}
    When $-2<\s\le 0$, assume that $\int_{(\R^\d)^2}|\g|(x-y)d|\mu|^{\otimes 2}<\infty$ to ensure that the modulated energy is finite.
    Let $v$ be a Lipschitz vector field. Then for any pairwise distinct point configuration $\XN\in(\R^\d)^N$, the following hold.
    
        \cite{rosenzweig_sharp_nodate} If $\s\ge \max(0,\d-2)$, then
        \begin{align}\label{eq:FIsupC}
            | \As_1[\XN,\mu,v]|
            \le C\|\nab v\|_{L^\infty}\Big(\Fr_N(\XN,\mu) - \frac{\log \la}{2N}\indic_{\s=0} + C_p\|\mu\|_{L^p}\la^{\frac{\d(p-1)}{p}-\s} \Big).
        \end{align}

        \cite{hess-childs_sharp_2025} If $0\le \s<\d-2$, then for any $\as \in (\d,\d+2)$, 
        \begin{multline}\label{eq:FIsubC1}
            | \As_1[\XN,\mu,v]| \\
            \le             C(\|\nab v\|_{L^\infty} + \|\Dm^{\frac{\as}{2}}v\|_{L^{\frac{2\d}{\as-2}}}\indic_{\as>2}) \Big(\Fr_N(\XN,\mu) - \frac{\log \la}{2N}\indic_{\s=0} + C_p\|\mu\|_{L^p}\la^{\frac{\d(p-1)}{p}-\s} \Big),
        \end{multline}
        \cite{nguyen_mean-field_2022,rosenzweig_global--time_2023} and for any $\frac{\d}{\d-\s-1}<p\le \infty$, 
        \begin{multline}\label{eq:FIsubC2}
            | \As_1[\XN,\mu,v]|
           \\ \le C(\|\nab v\|_{L^\infty} + \|\Dm^{\frac{\d-\s}{2}}v\|_{L^{\frac{2\d}{\d-\s-2}}}\indic_{\s<\d-2} )\Big(\Fr_N(\XN,\mu) + C_p\|\mu\|_{L^p}\ka^{\frac{\d(p-1)}{p}-\s}(1-(\log\ka) \indic_{\s=0})\Big).
        \end{multline}
        
        \cite{rosenzweig_wasserstein_nodate} If $-2<\s <0$, then
        \begin{align}\label{eq:FInonsing}
             |\As_1[\XN,\mu,v]| \le C(\|\nab v\|_{L^\infty} + \|\Dm^{\frac{\d-\s}{2}}v\|_{L^{\frac{2\d}{\d-\s-2}}}\indic_{\s<\d-2})\Fr_N(\XN,\mu)
        \end{align}
    In all cases, the constant $C>0$ depends only on $\d,\s$ and the constants $C_\as, C_p>0$ additionally depend on $\as,p$, respectively.    
    \end{thm}

{
\begin{remark}
Here, optimality is understood in the ``worst case'' sense, i.e. that there exists a point configuration $\ux_N$ that achieves the bound. This should be contrasted with the ``average'' sense, where $\ux_N\sim f_N \in \mathcal{P}((\R^\d)^N)$ and one estimates $\E_{\ux_N\sim f_N}[|\eqref{15}|]$. It is well-known that averaging may produce smaller errors. For instance, if $f_N \sim \mu^{\otimes N}$ (i.e. $\mu$-iid), then $\E_{\ux_N\sim f_N}[|\eqref{15}|]\approx N^{-1}$. More generally, one can prove \emph{high-temperature/entropic} commutator estimates for bounded interactions \cite{jabin_quantitative_2018, LLN2020} and even some Riesz cases \cite{DGR}.  
\end{remark}
}
{
\begin{remark}
\cref{thm:FIprime} (and its proof) hold mutatis mutandis on the flat torus $\T^\d$ where $\g$ is now taken to be the solution of $\Dm^{\d-\s}\g = \cd(\delta_0-1)$. 
\end{remark} 
}

Given a measure $\nu$ with density in $\Sc(\R^\d)$ and $\int_{\R^\d}d\nu=1$, letting $f\coloneqq \nu-\mu$, we can find a sequence  of pairwise distinct point configurations $\XN$ such that, as $N \to \infty$,
    \begin{align}
        \Fr_N(\XN,\mu)\rightarrow \frac12\int_{(\R^\d)^2}\g(x-y)f(x)f(y)dxdy = \frac{\cd}{2}\|f\|_{\dot{H}^{\frac{\s-\d}{2}}}^2
    \end{align}
    and
    \begin{align}
        \As_1[\XN,\mu,v] \rightarrow \frac12\int_{(\R^\d)^2}\nab\g(x-y)\cdot (v(x)-v(y))f(x)f(y)dxdy =\frac12\Big\langle f,\comm{v\cdot}{\nab \Dm^{\s-\d}}f\Big\rangle_{L^2}.
    \end{align}
    Thus, a commutator estimate with an additive error involving the singular centered measure $\frac1N\sum_{i=1}^N\delta_{x_i}-\mu$ implies a commutator estimate without an additive error for the regular centered measure $f$. More precisely, \cref{thm:FIprime} implies the following  standard or \emph{unrenormalized} commutator estimate (see \cite[Proposition 3.1]{nguyen_mean-field_2022}).
    
    \begin{prop}\label{prop:FI'}
        Let $\d\geq 1$ and $-2<\s<\d$. Then there exists $C>0$, depending only on $\d$ and $\s$, such that for all $f,g\in\mathcal{S}(\R^\d)$,\footnote{If $\s=0$, then we imolicitly assume $f,g$ are zero-mean to avoid the low-frequency divergence in $\dot{H}^{-\frac{\d}{2}}$.} we have
		\begin{multline}\label{eq:CE}
			\bigg|\int_{(\R^\d)^2} (v(x)-v(y))\cdot\nabla\g(x-y)f(x)g(y)dxdy\bigg|
			\\\leq C\Big(\|\nabla v\|_{L^\infty}+\||\nabla|^{\frac{\d-\s}{2}} v\|_{L^\frac{2\d}{\d-\s-2}}\indic_{\s<\d-2}\Big)\|f\|_{\dot{H}^{\frac{\s-\d}{2}}}\|g\|_{\dot{H}^{\frac{\s-\d}{2}}}.
		\end{multline}
    \end{prop}

    
    In fact---as we describe below---the philosophy of \cite{rosenzweig_mean-field_2020, nguyen_mean-field_2022, rosenzweig_sharp_nodate}, and to a more implicit degree \cite{leble_fluctuations_2018,  duerinckx_mean-field_2016, serfaty_mean_2020, serfaty_gaussian_2023}, proceeds in the reverse order: first show the {unrenormalized} commutator estimate of \cref{prop:FI'}, then combine it with a renormalization procedure to allow for the irregular test measure $\frac1N\sum_{i=1}^N\delta_{x_i}$. It is this renormalization step, which is the most technically involved part of the proof, that leads to the additive $N$-dependent errors in the singular case $0\le \s<\d$. For this reason, the estimates \eqref{eq:FIsupC}, \eqref{eq:FIsubC1}, \eqref{eq:FIsubC2} are sometimes referred to as \emph{renormalized} commutator estimates to distinguish them from the unrenormalized estimate \eqref{eq:CE}. In the nonsingular case $-2<\s<0$, no renormalization is necessary, and one can apply a commutator estimate directly, leading to the absence of additive errors in \eqref{eq:FInonsing}.

We will only discuss the proofs of the (super-)Coulomb estimate \eqref{eq:FIsupC} (\cref{ssec:rcomsupC}) and the sub-Coulomb estimate  \eqref{eq:FIsubC1} (\cref{ssec:rcomsubC}), as these are two sharp, singular cases. The reader interested in the proof of the nonsharp sub-Coulomb estimate \eqref{eq:FIsubC2} may consult \cite{nguyen_mean-field_2022}. The estimate \eqref{eq:FInonsing}, shown in \cite{rosenzweig_wasserstein_nodate}, follows from a straightforward extension of the method of \cite{nguyen_mean-field_2022}.   

\subsection{Coulomb/super-Coulomb case}\label{ssec:rcomsupC}
Following the chronology of the literature, we start with the (super-)Coulomb case $\s\ge \max(0,\d-2)$ as considered in \cite{rosenzweig_sharp_nodate}.

As in prior works  \cite{duerinckx_mean-field_2016,leble_fluctuations_2018,serfaty_mean_2020}, the starting point is the {\it electric reformulation} of the modulated energy \eqref{hofis:eq:modenergy} as a (renormalized) version of the energy 
\begin{equation}\label{hofis:eq:gradHn}
\int_{\R^\d} |\nabla h_N|^2, \qquad h_N\coloneqq\g* ( \mu_N-\mu),
\end{equation}
where $h_N$ denotes the electric potential and as before, $\mu_N$ is the empirical measure. This reformulation is valid for the Coulomb case. In the super-Coulomb case, we may instead use the Caffarelli-Silvestre extension procedure  \cite{CS2007}. It allows to  view $\g$ as the kernel of a local operator $-\div (\zg \nab \cdot)$ and 
replace 
 \eqref{hofis:eq:gradHn}   by the  same quantity with weight $\zg$, once working in the extended space $\R^{\d+1}$, with $z$ being the $(\d+1)$-th variable.   The fact that such a procedure is not available in the same way  if $\s<\d-2$  (one has to consider Dirichlet-type energies with higher derivatives; see, e.g., \cite{yangHigherOrderExtensions2013,coraSpolyharmonicExtensionProblem2022}) is what restricts us  so far to the range $\s \ge \d-2$. 

The key to the proof of inequalities such as \eqref{eq:FIsupC} was the observation of a {\it stress-energy tensor structure} in the left-hand side. More precisely, letting
\begin{equation}\label{hofis:defTT}
(T_h)_{ij}\coloneqq  2\partial_i h \partial_j h - |\nabla h|^2 \delta_{ij}.
\end{equation}
denote the stress-energy tensor (from classical mechanics or calculus of variations) associated to $h$, it is well-known that if $h$ is smooth enough, then 
\begin{equation} \label{hofis:divt}
\div T_h= 2 \nabla h \Delta h,
\end{equation} 
where the {divergence may be taken with respect either rows or columns, as the tensor is symmetric.} Ignoring the question of the diagonal in \eqref{hofis:eq:modenergy}, the main point is that in the Coulomb case, using \eqref{hofis:divt}, one may rewrite the left-hand side of \eqref{eq:FIsupC} after desymmetrization as 
\begin{align}\label{hofis:sstruc} 
2\int_{\R^\d}\int_{\R^\d} v(x) \cdot \nabla \g(x-y) d(\mu_N- \mu)^{\otimes 2} &=  -\frac{2}{\cd}\int_{\R^\d} v\cdot \nabla h_N \Delta h_N =-\frac{1}{\cd} \int_{\R^\d} v \cdot \div T_{ h_N}.
\end{align}
An integration by parts then allows to formally conclude, since $\int_{\R^\d} |T_{ h_N}|\le \int_{\R^\d} |\nabla h_N|^2$, which is the energy. The preceding reasoning is purely formal because $|\nabla h_N|^2$ diverges near each  point of the configurations (this is related to the diagonal excision), and this computation needs to be  properly renormalized, which is the main technical roadblock in the proof. Things work the same in the Riesz case, after extending the space and adding the appropriate weight. We refer to \cite[Section 3]{rosenzweig_sharp_nodate} for the details.


Considering a pure stress-tensor approach to the proof of higher-order estimates \eqref{hofis:mainn} for $n \ge 2$, it is unclear that an algebraic manipulation like \eqref{hofis:sstruc} can be found, leading to delicate proofs in \cite{leble_fluctuations_2018,serfaty_gaussian_2023} which do not seem extendable to $n \ge 3$. In \cite{rosenzweig_sharp_nodate}, we are able to exhibit a suitable---albeit more complicated---{stress-tensor structure in the higher-order variations} \eqref{15}. This structure requires us to consider not just the electric potential $h_N$, but also {\it iterated commutators} of it. Let us describe this more precisely at first order in the Coulomb case. We will discuss the higher commutators later in \cref{ssec:rcomloc}.

Given a {measure} $f$ in $\R^\d$ of integral $0$ (think of $f=\mu_N-\mu$), $h^f = \g* f$ its Coulomb potential, and a vector-field $v$, we define the {\it first commutator} of $h^f$ as 
\begin{equation}
\kappa^{f}(x) \coloneqq \int_{\R^\d} \nabla \g(x-y) \cdot (v(x)-v(y)) df(y).
\end{equation}
It is a commutator  because it can be rewritten as 
\begin{equation}\label{hofis:rewrikappa}
\kappa^f =-\g * (\div (vf)) +v \cdot \nabla \g  * f= - h^{\div (vf)}+ v \cdot \nabla h^f.
\end{equation}
In fact, the commutator  $\kappa^f$ is intimately tied to the stress-tensor itself, via the relation \eqref{hofis:sstruc}, which can be rewritten in general as
\begin{equation}\label{hofis:divtk}
2\cd\int_{\R^\d} \kappa^f dw = \cd\int_{\R^\d} \kappa^f dw  + \cd\int_{\R^\d} \kappa^w df  =  -\int_{\R^\d} v \div \comm{\nab h^f}{\nab h^w}= \int \nabla v : \comm{\nab h^f}{\nab h^w},
\end{equation}
where
\begin{align}\label{hofis:eq:introst}
\comm{\nab h^f}{\nab h^w}_{ij} \coloneqq \p_i h^f\p_j h^w + \p_j h^f \p_i h^w - \nab h^f\cdot\nab h^w\delta_{ij}	
\end{align}
is the bilinear generalization of the stress-energy tensor $T_{h}$ from \eqref{hofis:defTT}, which we note is symmetric. 

Considered separately, each term in the definition of the function $\kappa^f$ is one derivative less regular than $h^f$; however thanks to its commutator structure, some compensation happens.  By choosing $w= -\frac{1}{\cd}\Delta \ka^{f}$, integrating by parts on the left-hand side, and using Cauchy-Schwarz on the right-hand side, we can
show {the $\dot{H}^1$ control}
\begin{equation}\label{hofis:commcontrol1}
\int_{\R^\d} |\nabla \kappa^f |^2  \le C \|\nabla v\|_{L^\infty} \int_{\supp \nab v} |\nabla h^f|^2.
\end{equation}
 Compared to classical proofs of commutator estimates via paraproducts, the proof of \eqref{hofis:commcontrol1} is remarkably simple, thanks to the use of the stress tensor. Moreover, the right-hand side is localized to the support of the transport, something not evidently obtainable with classical proofs. 

\subsection{sub-Coulomb case}\label{ssec:rcomsubC}
We now discuss the sub-Coulomb case $\s<\d-2$ following the later work \cite{hess-childs_sharp_2025} joint with Hess-Childs and Serfaty.  In fact, the cited work gives a new proof for the full Riesz case $0\le \s<\d$, but the estimates are worse in the (super-)Coulomb case in terms of their transport regularity assumptions. We comment more on this point in \cref{ssec:rcomtransreg}.

There are two main ingredients in the proof: (1) a new potential truncation scheme based on an integral representation of the Riesz potential and  (2) commutator estimates of Kato-Ponce type.

The starting point for the first ingredient is the identity
\begin{align}\label{eq:introgintrep}
\forall x\ne 0,\qquad \g(x) = \mathsf{c}_{\phi,\d,\s}\int_0^\infty t^{\d-\s}\phi_t(x)\frac{dt}{t},
\end{align}
where $\phi$ is any sufficiently nice radial function with nonnegative Fourier transform $\hat{\phi}\ge 0$, $\phi_t(x)\coloneqq t^{-\d}\phi(x/t)$, and $c_{\phi,\d,\s}$ depends on $\phi,\d,\s$. In the log case $\s=0$, the integral has to be renormalized at large scales to obtain a convergent expression due to the fact that $\log|x|$ does not decay at infinity. The reader may consult \cite[Section 2]{hess-childs_sharp_2025} for precise assumptions ensuring the validity of the identity. The formula, which is a consequence of scaling, expresses that the Riesz potential is a weighted average of the family of approximate identities $\{\phi_t\}_{t>0}$. In \cite[Chapter 6]{Rubin2024}, Rubin refers to this as a ``wavelet type representation.'' The identity \eqref{eq:introgintrep} may also be understood in terms of the scaling/inversion properties of Mellin transforms \cite[Appendix B]{FS2009}.

The utility of the representation \eqref{eq:introgintrep} is that it conveniently allows us to truncate the singularity of the Riesz potential simply by cutting off the integral for small $t$. More precisely, given a scale $\eta>0$, we define the \emph{truncated potential}
\begin{align}\label{eq:introgeta}
\g_\eta(x) \coloneqq  \mathsf{c}_{\phi,\d,\s}\int_\eta^\infty t^{\d-\s}\phi_t(x)\frac{dt}{t}.
\end{align}
The truncated potential has several important properties that are crucial to the analysis (see \cite[Lemma 2.1]{hess-childs_sharp_2025} for all of the properties of $\g_\eta$):
\begin{enumerate}
\item While $\g$ is singular at the origin, $\g_\eta$ is finite and bounded by $\g(\eta)$.
\item $\g_\eta$ is a positive definite kernel (i.e. repulsive, or with nonnegative Fourier transform).
\item $\g_\eta \le \g $ and the difference $\g-\g_\eta$ decays rapidly at scale $\eta$ outside the ball $B(0,\eta)$.
\end{enumerate}

Without using the representation \eqref{eq:introgintrep}, it is not easy to construct a truncation $\g_\eta$ satisfying all of the three requirements at once, and  
a standard mollification  procedure fails to do so. Moreover, the choice of $\phi$ will require care.
Thanks to these properties,  splitting $\Fr_N(\ux_N,\mu)$ into 
\begin{multline}
\Fr_N(\ux_N, \mu)= \frac12 \int_{(\R^\d)^2}\g_\eta(x-y)d\Big(\frac1N\sum_{i=1}^N \delta_{x_i}-\mu\Big)^{\otimes 2}\\
+ \frac12 \int_{(\R^\d)^2\backslash \triangle }(\g-\g_\eta)(x-y)d\Big(\frac1N\sum_{i=1}^N \delta_{x_i}-\mu\Big)^{\otimes 2} - \frac{\g_\eta(0)}{2N},
\end{multline}
where the diagonal can be reinserted in the first integral up to the well-controlled last term on the second line, we easily obtain two sets of controls by $\Fr_N$. The first is on the first term on the right-hand side, which is nonnegative (thanks to the repulsive nature of $\g_\eta$) and is the square of a genuine MMD. 
The second is a control on the second term on the right-hand side which, up to well-controlled additive error terms, is also nonnegative and, thanks to properties of $\g_\eta$, controls $\frac{1}{N^2} \sum_{i\neq j} \g(x_i-x_j)$, i.e.,~the interaction energy coming from particle interactions at small scales.  See \cite[Proposition 2.11]{hess-childs_sharp_2025} for the precise statements of these controls. The former control allows to show in  \cite[Proposition 2.15]{hess-childs_sharp_2025} that the modulated energy controls the squared inhomogeneous Sobolev norm $\|\frac1N\sum_{i=1}^N\delta_{x_i}-\mu\|_{H^{-\frac{\d}{2}-\varepsilon}}^2$, for any $\varepsilon>0$, up to $O(\la^{\d-\s})$ error (i.e. coercivity). Choosing $\eta \propto N^{-1/\d}$ (the microscale), the latter control allows to bound in terms of the modulated energy the interaction energy due to nearest neighbors at the microscale, a fortiori bounding the number of points with nearest-neighbor distance below the microscale, see \cite[Corollary 2.14]{hess-childs_sharp_2025}. Let us mention that the potential truncation introduced here plays a key role in our forthcoming work \cite{hess-childs_optimal_nodate} on optimal quantization via MMDs with kernels having power law-type spectra (e.g. Riesz).



In contrast to \cite{petrache_next_2017, serfaty_mean_2020, nguyen_mean-field_2022, rosenzweig_sharp_nodate}, there is no need to consider a smearing/regularizing of the charges $\delta_{x_i}$, nor to use smearing radii that are  point-dependent and configuration-dependent---we only need to consider the truncated potential $\g_\eta$.
In these previous works, the strategy behind an estimate of the form in \cref{thm:FIprime} is to first prove an estimate valid when the diagonal is reinserted and $\frac{1}{N}\sum_{i=1}^N\delta_{x_i}-\mu$ is replaced by a more regular distribution $f$. One then reduces to this regular setting by replacing the point charges $\delta_{x_i}$ by smeared charges $\delta_{x_i}^{(\eta_i)}$ and estimating the error from this replacement, the so-called renormalization step. Important to implementing the renormalization is to show that the modulated energy controls both the potential energy of the difference $\frac1N\sum_{i=1}^N\delta_{x_i}^{(\eta_i)}-\mu$ and the small-scale interactions.

The strategy in \cite{hess-childs_sharp_2025} differs crucially in two regards. First, rather than regularize $\frac{1}{N}\sum_{i=1}^N\delta_{x_i}-\mu$ through smearing, we regularize $\g$ through the truncation $\g_\eta$.\footnote{The work \cite{bresch_modulated_2019} also directly regularized the interaction $\g$ at small length scales. However, the procedure in that work is completely different than in \cite{hess-childs_sharp_2025}, limited to the spatially periodic setting, and does not yield optimal error estimates.} Second, to prove an estimate for  $\g_\eta$, we use the definition \eqref{eq:introgeta} to reduce to proving estimates for
\begin{align}
\int_{(\R^\d)^2} (v(x)-v(y))\cdot \nabla \phi\Big(\frac{x-y}{t}\Big)d\Big(\frac{1}{N}\sum_{i=1}^N\delta_{x_i}-\mu\Big)^{\otimes 2}(x,y) , \qquad t \in (0,\infty),
\end{align}
which we integrate over $(\eta,\infty)$ with respect to the measure $t^{\d-\s}\frac{dt}{t}$. Provided the dependence on $v$ in our estimate is scale-invariant, this reduces to proving an estimate for $t=1$. 

We now come to the second main ingredient of our proof. The identity \eqref{eq:introgintrep} is valid for a large class of scaling functions $\phi$, but not all choices of $\phi$ will yield the desired estimate. Indeed, it is a elementary Fourier computation that no estimate of the form
\begin{align}
\int_{(\R^\d)^2}(v(x)-v(y))\cdot\nab\phi(x-y)f(x)f(y) \le C_v \int_{(\R^\d)^2}\phi(x-y)f(x)f(y)
\end{align}
can hold when $\hat\phi$ decays super-polynomially (e.g. Gaussian). It turns out that  a good choice for $\phi$ is a Bessel potential, i.e. $\hat\phi(\xi) = \jp{2\pi\xi}^{-\as}$, which is the kernel of the inhomogeneous Fourier multiplier $\jp{\nab}^{\as}$, where we use the Japanese bracket notation $\jp{x}^{a} \coloneqq (1+|x|^2)^{a/2}$ for $a\in\R$ and $x\in\R^\d$. The Bessel potential is a screening of the Riesz potential which preserves its local behavior at the origin but avoids the issues at low frequency that lead to the slow spatial decay of the Riesz potential. Setting $h\coloneqq \phi\ast f$ and using that $\jp{\nab}^{\as}h = f$, the  choice of Bessel potential reduces to proving a product rule for $\jp{\nab}^{\as/2}$. Such product rules are known as Kato-Ponce estimates in the harmonic analysis literature on account of their origins \cite{KP1988}. This is the content of \cite[Proposition 3.1]{hess-childs_sharp_2025}, which is the main technical lemma. Here, we rely on commutator estimates of Li \cite{Li2019}, as well as a local representation of fractional powers $\jp{\nab}^{\as/2}$ via dimension extension (see \cite[Section 3]{hess-childs_sharp_2025}), analogous to the representation for the fractional Laplacian popularized by Caffarelli and Silvestre \cite{CS2007}.

The parameter $\as$ for the Bessel potential is subject to constraints. We need $\as>\d$, so that $\phi$ is continuous, bounded. On the other hand, if $\as$ is too large, then we will not be able to obtain an estimate whose dependence on $v$ scales properly. This leads to the constraint $\as<\d+2$. Even with the constraints on $\as$, it is somewhat remarkable that an intermediate inhomogeneous estimate allows for a homogeneous estimate in the end.

The identity \eqref{eq:introgintrep} itself suggests a class of potentials beyond the exact Riesz case by replacing the function $t^{\d-\s}$ with a general function of $t$ that is pointwise controlled above and below by $t^{\d-\s}$. We call these \emph{Riesz-type} potentials. The reader may find their proper definition in \cite[Section 2.2]{hess-childs_sharp_2025} and should compare this class with the related, but different, class of Riesz-type potentials previously introduced in \cite{nguyen_mean-field_2022}. In fact, the estimate \cref{thm:FIprime}\eqref{eq:FIsubC1} is obtained as a special case of a more general functional inequality for Riesz-type potentials (see \cite[Theorem 4.1]{hess-childs_sharp_2025}).

\subsection{Localized and higher-order estimates}\label{ssec:rcomloc}
 Having discussed first-order estimates, let us now see how to obtain the second- and higher-order, possibly localized, estimates using stress-tensor and commutator structures, following \cite{rosenzweig_sharp_nodate}. The results discussed in this subsection are limited to the (super-)Coulomb case. It is unclear how to obtain (sharp) higher-order estimates in the sub-Coulomb case with the method of \cite{hess-childs_sharp_2025}. For more discussion on this problem, see \cite[Section 1.5]{hess-childs_sharp_2025} or \cite[Sections 5.2, 5.3]{hess-childs_sharp_2025}.
 
 Ignoring the question of renormalization, and sticking to the Coulomb case, we have seen in \eqref{hofis:sstruc} (combining with \eqref{15}) that the first variation of the modulated energy along the transport map $\I+tv $ is
\begin{align}\label{hofis:124}
-\frac1{2\cd}\int_{\R^\d} v \cdot \div T_{h_N}dx = \frac1{2\cd} \int_{\R^\d} \nabla v : T_{h_N}dx.
\end{align}
To compute the second variation of the modulated energy, we thus need to compute the first variation of $\div T_{h_N}$ when again $\mu_N$ and $\mu$ are pushed forward by $\I+ tv$. In view of the expression for $T_{h_N}$ in \eqref{hofis:defTT}, it suffices to compute the derivative of $h_N^t$ (with obvious notation) at $t=0$, and since $
h_N= \g * (\mu_N-\mu)$, the definition of the push-forward yields that 
\begin{equation} \label{hofis:dthnt}\frac{d}{dt}\big|_{t=0} h_N^t= -\g * (\div (v (\mu_N-\mu)) =-  h^{\div (v(\mu_N-\mu))}=  \kappa^{\mu_N-\mu} -v\cdot \nab h_N,\end{equation} after using \eqref{hofis:rewrikappa}. 
{The $L^2$ norm of the gradient of this expression is one derivative more singular than the energy $\int_{\R^\d} |\nabla h_N|^2$. }
Still with $f= \mu_N-\mu$, inserting \eqref{hofis:dthnt} into the variation of \eqref{hofis:124}, we can  decompose the second-order variation as 
\begin{equation}\label{hofis:line12}
\frac1{\cd}\int_{\R^\d} \p_i v^j : \left( \comm{\nab  h_N }{\nab  \kappa^f }-  \comm{\nab  h_N }{ \nab ( v\cdot \nabla h_N  ) }\right).
\end{equation}
Thanks to the $\dot{H}^1$ estimate \eqref{hofis:commcontrol1} for the commutator $\ka^f$, the  first  term on the right-hand side can be directly controlled by $ C_v \int |\nabla h_N|^2$ as desired, while the second one can be transformed into  similarly controllable terms---albeit $C_v$ is now quadratic in $v$ and depends on $v,\nabla^{\otimes 2}v$---by means of integration by parts of the advection operator $v \cdot \nabla $.

This argument can then be iterated at next order by introducing the family of $n$-th order  commutators $\kappa^{(0), f}=h^f$ and 
\begin{equation}
\kappa^{(n),f}\coloneqq \int_{\R^\d} \nabla^{\otimes n} \g(x- y) : (v(x)-v(y))^{\otimes n} df(y),
\end{equation}
together with ``transported" commutators
\begin{equation}
\kappa_t^{(n),f}\coloneqq \int_{\R^\d} \nabla^{\otimes n} \g(x-y-tv(y)) : (v(x)-v(y))^{\otimes n} df(y),
\end{equation}
and observing the elegant recursive formula
\begin{equation}\label{hofis:iteratek}
\partial_t \kappa_t^{(n),f} = \kappa_t^{(n+1),f} - v \cdot \nabla \kappa^{(n),f}_t,
\end{equation}
where $t$ is the same parameter as in the transport map $\I+tv$. In other words, the high-order ``time-dependent'' commutators $\ka_t^{(n),f}$ satisfy a hierarchy of transport equations with a source coupling the $(n+1)$-th order commutator to the $n$-th order commutator.  An iteration of the same argument  as for \eqref{hofis:commcontrol1} using this recursion allows to prove the estimate 
\begin{equation}\label{hofis:commuordern}
\int_{\R^\d} |\nabla \kappa^{(n),f}|^2 \le C_v \int_{\supp\nab v} |\nabla h^f|^2,
\end{equation}
valid for general functions $f$ such that the right-hand side is finite, with the constant $C_v $ now being $n$-linear in $v$ and involving derivatives up to order $n$ of $v$.  One can see the relation \eqref{hofis:commuordern} as an {\it $L^2$-based  regularity theory for arbitrary-order commutators}. 


The discussion so far has been limited to the Coulomb case,  the complete details of which are presented in \cite[Appendix A]{rosenzweig_sharp_nodate}. In generalizing to the Riesz case $\s \neq \d-1,\d-2$, we run into issues. The vector field $v$ may be trivially extended to $\R^{\d+\k}$, with $\k=1$ if $\s>\d-2$, by fixing the last component to zero, $\g$ may be trivially extended through radial symmetry, and all distributions $f$ on $\R^{\d}$ viewed as living on the boundary $\R^{\d}\times \{0\}^\k$. Consequently, $\ka^{(n),f}$ may be viewed as a function on $\R^{\d+\k}$.  Going through the computations above, all integrals $\int_{\R^\d}$ should be replaced $\int_{\R^{\d+\k}}|z|^{\gamma}$ where $\gamma=\s+2-\d-\k$. The obstruction lies in obtaining the $L^2$ gradient bound \eqref{hofis:commuordern} by duality.

Setting $L \coloneqq -\div(\zg\nab\cdot)$, we would like to take $w = \frac{1}{\cd} L\ka^{(n),f}$ in  $\int_{\R^{\d+\k}}\ka^{(n),f}w$
and then integrate by parts to conclude an estimate for $\int_{\R^{\d+\k}}\zg |\nab \ka^{(n),f}|^2$. However, this choice for $w$ is {a priori} not supported on $\R^{\d}\times\{0\}^\k$. If $\supp \phi \not\subset \R^{\d}\times \{0\}^\k$, it is not necessarily true that $\g\ast L\phi = \cd\phi$, for a test function $\phi$ on $\R^{\d+\k}$, unless $\ga = 0$. This forces us to come up with a new approach.

The starting point is the observation $L\g = \cd\delta_0$ in $\R^{\d+\k}$, which underlies the Caffarelli-Silvestre  representation of the fractional Laplacian as a degenerate elliptic operator in $\R^{\d+\k}$. From this observation, it follows that $\ka^{(n),f}$ satisfies the equation
\begin{multline}\label{hofis:eq:introLkanf}
L\ka^{(n),f} = \cd(-1)^n\sum_{\sigma \in \Ss_n} \p_{i_{\sigma_1}}{v}^{i_1}\cdots\p_{i_{\sigma_n}}v^{i_n}f\delta_{\R^\d\times\{0\}^\k}   -n\p_i(\zg\p_i v \cdot\nu^{(n-1),f})  \\
- n\zg\p_i v \cdot \p_i\nu^{(n-1),f}+ n(n-1)\zg(\p_i v)^{\otimes 2}:\mu^{(n-2),f},
\end{multline}
where $\mathbb{S}_n$ denotes the symmetric group on $[n]$ and a repeated index denotes summation over that index. In the right-hand side, $\nu$ is a vector field on $\R^{\d+\k}$ that ``morally'' is like $\nab\ka^{(n),f}$, while $\mu$ is a symmetric matrix field on $\R^{\d+\k}$ that morally is like $\nab^{\otimes 2}\ka^{(n),f}$ (see \cite[(4.12), (4.14)]{rosenzweig_sharp_nodate}  respectively, for the precise definitions).  The right-hand side of \eqref{hofis:eq:introLkanf} is good because it only depends on lower-order (i.e.~ $n-1$, $n-2$) commutators. Following the standard method for the $L^2$ regularity of elliptic equations, we want to obtain an estimate for $\int_{\R^{\d+\k}}\zg|\nab\ka^{(n),f}|^2$ by testing the equation \eqref{hofis:eq:introLkanf} against $\ka^{(n),f}$ and using an induction hypothesis that will guarantee good bounds on the lower order commutators and their derivatives. There is however a difficulty involving the nondivergence form of the right-hand side of \eqref{hofis:eq:introLkanf} and $L^2$ estimates for $\ka^{(n),f}$ (see the beginning of \cite[Section 4.1]{rosenzweig_sharp_nodate} for elaboration),  which requires the use of a Poincar\'e inequality when $n \ge 3$ and yields an estimate that deteriorates when $\ell$, the size of the support of $v$, gets large. 

As a result, when $\supp v$ is not compact, which is not the case for the aforementioned CLT application, this  approach does  not work, and the right-hand side of \eqref{hofis:mainn} is infinite. However, obtaining unlocalized $L^2$ commutator estimates is much easier and may be done following  the method of \cite{nguyen_mean-field_2022}: integration by parts and the difference quotient characterization of the Sobolev seminorm $\dot{H}^{\frac{\s-\d}{2}}$ when $\s>\d-2$, and the $L^2$ boundedness of Calder\'{o}n $\d$-commutators when $\s=\d-2$. 

For the precise statements of the results described above, we refer the reader to \cite[Theorem 4.1]{rosenzweig_sharp_nodate} (and more generally, Section 4.1), which is our omnibus $L^2$ regularity estimate for commutators. 

This discussion has so far left aside the question of renormalization, which is that of dealing with the singularities of the Dirac masses in $h_N$, and which is arguably the most delicate part of the analysis. As mentioned in \cref{ssec:rcomsupC} for the first-order estimates, the renormalization is handled as in \cite{serfaty_mean_2020} and subsequent works via a charge smearing and potential truncation at a lengthscale crucially {\it depending on each point} $x_i$ and proportional to its nearest-neighbor distance. 
The renormalization in the $n\ge 2$ case is much more delicate, occupying the largest part of the paper and requiring several innovations.

First, replacing the Diracs in the left-hand side of \eqref{hofis:mainn} by their smearings $\delta_{x_i}^{(\eta_i)}$ will not allow us to apply our unrenormalized commutator estimate unless $\s=\d-2$ because $\mu\delta_{\R^\d\times\{0\}^\k} - \frac1N\sum_{i=1}^N \delta_{x_i}^{(\eta_i)}$ is not supported on $\R^\d\times \{0\}^\k$. To overcome this difficulty, we introduce a new smearing $\rho_{x_i}^{(\eta_i)}$, obtained simply by mollifying the Dirac in $\R^\d$, and estimate the electric energy difference of the two regularizations in terms of the modulated energy. We can then directly apply our $L^2$ commutator bound to $\mu-\frac1N\sum_{i=1}^N\rho_{x_i}^{(\eta_i)}$. A novelty compared to previous work is that we will choose the mollifier to have vanishing moments up to high enough order in order to obtain optimal error rates. For the details of the described construction and estimates, see \cite[Section 5.1]{rosenzweig_sharp_nodate}.

Second, to handle the error from this smearing, we establish an {\it $L^\infty$-based regularity theory for commutators} (see \cite[Sections 4.5-4.6]{rosenzweig_sharp_nodate}), which shows that the $L^\infty$ norm of ``horizontal'' derivatives $\nab_{x}^{\otimes m}\kappa^{(n),f}$, where $\nab_x \coloneqq (\p_1,\ldots,\p_\d)$, {\it restricted to arbitrarily small balls} can be estimated in terms of the $L^2$ norm of $\nab\ka^{(n),f}$ in a double ball. Commutator estimates do not typically allow such localization, and therefore this result may be of independent interest. The proof relies on regularity theory for second-order divergence-form elliptic operators with $A_2$ weights \`{a} la Fabes et al. \cite{FKS1982} in the extended space $\R^{\d+\k}$, which, by the Caffarelli-Silvestre extension, may be viewed as a regularity theory for the fractional Laplacian. The proof also crucially uses the commutator structure in the form of equation \eqref{hofis:eq:introLkanf} for $L\ka^{(n),f}$, along with recursions for $\ka^{(n),f}$, $\nu^{(n),f}$, and $\mu^{(n),f}$. When $\ell \asymp 1$, this local regularity theory may be avoided with cruder estimates, which have the benefit of requiring less regularity on the vector field $v$. However, when $\ell\ll 1$, these cruder estimates give errors that scale sub-optimally.

Third, we improve truncation errors estimates from \cite{serfaty_mean_2020}, obtaining the optimal dependence on $\mu$, and show new localizable {\it mesoscale interaction energy estimates}, which generalize those of  \cite{serfaty_gaussian_2023} to the Riesz case. See \cite[Sections 2.2-2.3]{rosenzweig_sharp_nodate}.

The complete details of the renormalization step may be found in \cite[Section 5.2]{rosenzweig_sharp_nodate}  and more generally, Section 5 of the cited work. Having discussed the proof strategy and main innovations, we conclude this subsection with statements of the main results. 

The first theorem is a global estimate, which does not require any localization of the vector field $v$. Note that the additive error term is  $\lambda^{\d-\s} \propto (N^{-1+\frac{\s}{\d}})$,  which is the announced optimal estimate. Compared to the localized estimate \eqref{hofis:mainn} below, the regularity requirement for $v$ is not dependent on the order $n$ of the commutator. As the reader may check, if $\s>\d-1$, we only need $\|\nab v\|_{L^\infty}$ ($m=0$) to achieve the sharp error $\la^{\d-\s}$. If $\d-2<\s\le \d-1$, then we also need $\|\nab^{\otimes 2}v\|_{L^\infty}$ ($m=1$). If $\s=\d-2$, then we additionally need $\|\nab^{\otimes 3}v\|_{L^\infty}$ ($m=2$).

{
\begin{thm}[Global estimates]\label{hofis:thm:mainunloc}
There exists a constant $C>0$ depending only $\d,\s$ such that the following holds. Let $\mu \in L^1(\R^\d)\cap L^\infty(\R^\d)$ with $\int_{\R^\d}\mu=1$. {If $\s\le 0$, suppose further that $\int_{(\R^\d)^2}|\g(x-y)|d|\mu|^{\otimes 2}<\infty$.} Let  $v:\R^\d\rightarrow\R^\d$ be  a Lipschitz vector field, suppose that $\la<1$, where $\lambda:= (N\|\mu\|_{L^\infty(\R^\d)})^{-\frac1\d}$.  
For any integers $m\ge 0$ and $n\ge 1$, there exists a constant $C>0$ depending only $\d,\s,m,n$ such that the following holds. Let  $v:\R^\d\rightarrow\R^\d$ be a ${C}^{m,1}$ vector field.  For any pairwise distinct configuration $\ux_N \in (\R^\d)^N$, we have
\begin{multline}\label{hofis:eq:HOFIunloc}
|\As_n[X_N,\mu,v]|  \le C(\|\nab v\|_{L^\infty})^n\Big(\Fr_N (\ux_N,\mu) -  { \Big(\frac{\log \lambda}{2N}\Big) }\indic_{\s=0} +  C\|\mu\|_{L^\infty}\lambda^{\d-\s}\Big)\\
 + C\sum_{q=1}^{m+1}\sum_{r=0}^{{ q}}D_{v,q,r}\Bigg({\la^{m+1}} +\la^{m+1-q}\Big(\Fr_N (\ux_N,\mu) -  { \Big(\frac{\log \lambda}{2N}\Big) }\indic_{\s=0} +  C\|\mu\|_{L^\infty}\lambda^{\d-\s}\Big) \\
 + {\la^{m+1}\|\mu\|_{L^\infty} } \begin{cases} 1, &  \s+q\neq \d \\
\log (1/\lambda),  & \s+q=\d \end{cases} 
+{\la^{m+1}}\begin{cases}\|\mu\|_{L^\infty}\la^{\d-\s-q}, & \s+q>\d \\ \|\mu\|_{L^\infty}(1 - \la^{\d-\s-q}) +  \|\mu\|_{L^1}, & {\s+q<\d} \\ \|\mu\|_{L^\infty}\log(1/\la) + \|\mu\|_{L^1}, & {\s+q=\d} \end{cases}\Bigg),
\end{multline}
where
\begin{align}
D_{v,q,r} \coloneqq  \sum_{\substack{1\le c_1,\ldots,c_r \\ c_1+\cdots+c_r = m+1-(q-r)}}\|\nab^{\otimes c_1}v\|_{L^\infty}\cdots\|\nab^{\otimes c_r}v\|_{L^\infty}\|\nab v\|_{L^\infty}^{n-r}.
\end{align}
\end{thm}
}

\smallskip
In studying the fluctuations of Coulomb/Riesz gases on mesoscopic lengths $N^{-1/\d}\ll \ell \ll 1$, it is crucial that commutator estimates admit a \emph{local} form, in which the right-hand side can be expressed in terms of the modulated energy restricted to the region where $v$ is supported (up to boundary errors). Such estimates are key analytic input to the transport method for proving CLTs \cite{leble_fluctuations_2018,serfaty_gaussian_2023,peilen_local_2025} (see \cite[Part 3]{serfaty_lectures_2024} for a consolidated treatment).

Our second theorem is a localized generalization of \cref{thm:FIprime}, as well as a higher-order version. Now, the modulated energy on the right-hand side is localized to a region comparable to the support of $v$ and the estimate scales with the diameter of the support. The additive error term is in $N^{-1}\lambda^{\d-\s} \propto (N^{-2+\frac{\s}{\d}})$ per point,  which is optimal. However, the error now depends on higher-order derivatives of $v$ and deteriorates when the size of the support of $v$ gets large. We mention that \cref{hofis:thm:FI} was an important component to the recent work \cite{peilen_local_2025} of Peilen and Serfaty on the fluctuations of super-Coulomb Riesz gases. 

To formulate the result, let $\Omega \subset\R^\d$, meant to represent $\supp v$. Associated to $\Omega$, define the microscopic length scale 
	\begin{equation}\label{hofis:deflambda}
	\lambda \coloneqq (N\|\mu\|_{L^\infty(\Omega)})^{-\frac{1}{\d}},
	\end{equation}
 which can be thought of as the typical inter-particle distance,
and define $\hat \Omega := \{x:\dist(x,\Omega)\le \frac14\la\}$. 
We will use $\ell$ to denote the typical size of $\Omega$ (or of $\supp v$), which may also depend on $N$ and whose only constraint is  to remain larger than the microscopic scale $\lambda$.
Finally, we let {$I_\Omega\coloneqq \{1\le i\le N: x_i\in\Omega\}$} and use $\#$ to denote the cardinality of a finite set.



\begin{thm}[Localized estimates]\label{hofis:thm:FI}
There exists a constant $C>0$ depending only $\d,\s$ such that the following holds. Let $\mu \in L^1(\R^\d)\cap L^\infty(\R^\d)$ with $\int_{\R^\d}\mu=1$. {If $\s\le 0$, suppose further that $\int_{(\R^\d)^2}|\g(x-y)|d|\mu|^{\otimes 2}<\infty$.} Let  $v:\R^\d\rightarrow\R^\d$ be  a Lipschitz vector field  and $\Omega$ be a closed set  containing a $2\lambda$-neighborhood of $\supp \nab v$, where $\lambda$ is  defined as in \eqref{hofis:deflambda} and $\la<1$.  For any pairwise distinct configuration $\ux_N \in (\R^\d)^N$, it holds that
\begin{align}\label{hofis:main1}
|\As_1[X_N,\mu,v]| \leq C\|\nabla v\|_{L^\infty}\Big({\int_{\Omega\times [-\ell, \ell]^{\k}}\zg|\nab h_{N,\vec{\rs}}|^2} + C\frac{\# I_\Omega\|\mu\|_{L^\infty(\hat \Omega)}\la^{\d-\s}}{N} \Big).
\end{align}
{Here, $h_{N,\vec{\rs}}$ is the truncated electric potential and $\vec{\rs}=(\rs_i)_{i=1}^N$ is the nearest-neighbor type distance defined by
\begin{align}
h_{N,\vec{\rs}}\coloneqq \frac{1}{N}\sum_{i=1}^N \G_{\rs_i} (\cdot-x_i) - \G\ast \mu, \qquad \g_{\rs_i}\coloneqq \min(\g,\g(\rs_i)), \qquad
\rs_i \coloneqq \frac14\min\Big(\min_{j\neq i}|x_i-x_j|, \lambda\Big),
\end{align}
where we use the radial symmetry of $\g$ to abuse notation.
}

Suppose in addition  that  $v\in C^{2n-1,1}$ and that  $\Om'$ is a ball of radius $\ell$ containing a $\lambda$-neighboorhood of $\supp v$ and $\Omega$  contains a $5\ell$-neighborhood of $\Omega'$,\footnote{{For integer $k\ge 0$ and real $\al\in [0,1]$, $C^{k,\alpha}$ denotes the inhomogeneous H\"{o}lder space of $k$-times continuously differentiable functions whose $k$-th derivative is $\alpha$-H\"older continuous.}} where $\la,\ell$ satisfy $\min(\la,\la/\ell) < \frac{1}{2}$. For any $n \ge 2$, we have 
\begin{multline} \label{hofis:mainn}
|\As_n[X_N,\mu,v]| \le  C\sum_{p=0}^n (\ell\|\nab^{\otimes 2}v\|_{L^\infty})^p
 \sum_{\substack{1\leq c_1,\ldots,c_{n-p} \\ n-p\le c_1+\cdots+c_{n-p} \le 2n}} \lambda^{-(n-p)+\sum_{k=1}^{p} c_{n-k} } \\
 \times\|\nabla^{\otimes c_1} v\|_{L^\infty}\cdots\|\nabla^{\otimes c_{n-p}} v\|_{L^\infty} \Big({\int_{\Omega\times [-\ell,\ell]^{\k}} \zg |\nab h_{N,\vec{\rs}}|^2 }+  C\frac{\# I_\Omega\|\mu\|_{L^\infty(\hat \Omega)}\la^{\d-\s}}{N} \Big),
 \end{multline} 
where $C>0$ may depend on $n$ and the summation of the $c_i$ is understood as vacuous when $p=n$. 
\end{thm}

Note that the estimate \eqref{hofis:mainn} does not allow to take $\Omega=\R^\d$ because that would mean $\ell=\infty$, in which case the right-hand side is infinite. This is an unfortunate limitation of the result, one we believe is technical. We refer the interested reader to \cite[Section 4.1]{rosenzweig_sharp_nodate} and \cite[Section 5.3]{hess-childs_sharp_2025} for a discussion of the issue and possible ways to circumvent it.

\begin{remark}
While the conditions on $\Om,\Om',\la$ may appear circular upon first read, it is indeed possible for them to hold. See \cite[Remark 1.4]{rosenzweig_sharp_nodate} for details. 
\end{remark}

\subsection{Transport regularity and defective estimates}\label{ssec:rcomtransreg}
Besides the sharpness of the additive errors in these commutator estimates, discussed above, another important question, which we recently considered in the joint work \cite{hess-childs_sharp_2025} with Hess-Childs and Serfaty, is the optimal regularity assumption for the transport velocity $v$ in the aforementioned commutator estimates. Compared to the former, the latter question has received little attention in the literature.
    
It is natural to ask whether the Lipschitz assumption on $v$ may be weakened while still preserving the scaling. For example, can $\|\nab v\|_{L^\infty}$  be replaced by $\|\nab v\|_{BMO}$? It is also unclear whether, in the sub-Coulomb regime $\s<\d-2$, the additional term of the form $\|\Dm^{\frac{\as}{2}}v\|_{L^{\frac{2\d}{\as-2}}}$ is necessary and how small the exponent $\as$ may be taken. By Sobolev embedding, $\|\Dm^{\frac{\as}{2}}v\|_{L^{\frac{2\d}{\as-2}}} \lesssim \|\Dm^{\frac{\as'}{2}}v\|_{L^{\frac{2\d}{\as'-2}}}$ for $2(\d+1)>\as'\ge \as$. Moreover, there is a gap between our two sub-Coulomb estimates \eqref{eq:FIsubC1}, \eqref{eq:FIsubC2}. 

This regularity question is motivated by applications. It is a general principle in PDE that if an equation has a scaling invariance, then a natural function space to study the equation is one that is invariant under the equation's scaling \cite{klainerman_pde_2010}. Such \emph{scaling-critical} or \emph{scaling-invariant} spaces are expected to be thresholds for the well-posedness/ill-posedness of the equation along with other phenomena, such as finite-time blowup vs. global existence.

When $\mathsf{V}^t = 0$, the reader may check that if $\mu$ is a solution to the mean-field equation \eqref{eq:MFlim} with $\be=\infty$, then  $\tl{\mu}(t,x) \coloneqq  a\mu(bt,cx)$, for $a,b,c>0$, is also a solution provided that $bc^{\d-2-\s}=a$. Fixing $b=1$ so that the solution lifespan is preserved under rescaling, yields the relation $a=c^{\d-2-\s}$. 
For $1\leq p\leq\infty$, the reader may check that
\begin{equation}
\|\tl{\mu}\|_{L^p} = c^{\d-\s-2-\frac{\d}{p}} \|\mu\|_{L^p},
\end{equation}
which implies that $p=\frac{\d}{\d-\s-2}$ is the scaling-critical exponent. 
Thus, the potential can be at  most Coulomb in order for there to be a critical $L^p$ space. More generally, letting $\dot{W}^{a,p}$ denote the homogeneous Bessel potential (fractional Sobolev) space (see, e.g., \cite[Section 3.3]{stein_singular_1970} or \cite[Section 1.3]{grafakos_modern_2014}), we see that
$\dot{W}^{a,p}$ is a critical space if and only if $a=\frac{\d}{p} + \s+2-\d$. Other scales of function spaces, such as H\"older-Zygmund spaces (see, e.g., \cite[Section 4]{stein_singular_1970}), are also interesting to consider as critical spaces. 


When $\mu \in \dot{W}^{\frac{\d}{p} + \s+2-\d,p}$, for $1<p<\infty$, some basic Fourier analysis shows that the associated mean-field velocity $v=\M\nabla\g\ast\mu$ is in $\dot{W}^{\frac{\d}{p}+1,p}$.  Unfortunately, $\dot{W}^{\frac{\d}{p}+1,p}$ does not embed into the homogeneous Lipschitz space $\dot{W}^{1,\infty}$. If $f\in \dot{W}^{\frac{\d}{p}+1,p}$, then one only has $\|\nab f\|_{BMO} \leq C\|f\|_{\dot{W}^{\frac{\d}{p}+1,p}}$. Even when $\s+2-\d=k$ is an integer and taking $p=\infty$, having $\mu\in\dot{W}^{k,\infty}$ does not imply that $v\in \dot{W}^{1,\infty}$. This is problematic because all of the commutator estimates of \cref{thm:FIprime} feature $\|\nab v\|_{L^\infty}$ on the right-hand side. Thus, to show mean-field convergence when the mean-field density $\mu^t$ belongs to the critical space, it would be desirable to have a commutator estimate with a regularity assumption that is satisfied when $\mu\in \dot{W}^{\frac{\d}{p} + \s+2-\d,p}$.
	
    
In \cite{hess-childs_sharp_2025}, we showed that estimate~\eqref{eq:CE} can fail if either of the vector field norms is weakened. We also showed necessary conditions for analogous commutator estimates to hold for a broader class of interaction potentials. In particular, the norm controlling $v$ must involve enough derivatives to compensate for the high-frequency decay of $\hat\g$. This demonstrates a more general principle of regularity trade-off between the interaction and transport for commutator estimates. All of these results, which may be found in \cite[Section ~3]{hess-childs_sharp_2025}, are obtained through explicit counterexamples. Given {that} the renormalized commutator estimates of \cref{thm:FIprime} imply unrenormalized commutator estimates, the problem reduces to showing counterexamples when $\frac1N\sum_{i=1}^N\delta_{x_i}-\mu$ is replaced by a zero-mean Schwartz function.

The second contribution of \cite{hess-childs_sharp_2025}, motivated by applications to mean-field convergence at scaling-critical regularity, concerns substitutes for the functional inequality of \cref{thm:FIprime} when considering the mean-field vector field $ v =\M\nab\g\ast\mu$ for $\mu$ in a critical space. The results of the previous paragraph give a negative answer to the possibility of a commutator estimate holding in this case, regardless of the size of the additive error. That said, a \emph{defective} commutator estimate still holds, by which we mean that there exists a locally integrable function $\rho:[0,\infty)\rightarrow [0,\infty)$, called the \emph{defect}, such that for all $\ep>0$,
\begin{align}\label{eq:commdefec}
    |\As_1[X_N,\mu,v]| \le \rho(\ep)C_v\Big(\Fr_N(\XN,\mu) + C_\mu N^{-\alpha} \Big) + {\ep}^{\beta}N^{\gamma}C_v'\Big(\Fr_N(\XN,\mu) + C_\mu' N^{-\alpha} \Big),
\end{align}
where $\alpha > 0$, $\beta,\gamma\ge 0$ and $C_v,C_v',C_\mu,C_\mu'>0$ are constants depending on $\d,\s$ and quantitatively on norms of $v,\mu$, respectively. As we will explain later in \cref{ssec:MFreg}, if $\rho({\ep}) = O(|\log {\ep}|)$ as ${\ep}\rightarrow 0^+$, then a defective commutator estimate suffices to conclude an estimate for the modulated energy that in turn implies mean-field convergence. 

In \cite[Theorem 2.6]{hess-childs_sharp_2025}, we showed a new defective estimate  with $\rho(\ep) = O(|\log\ep|^{\theta})$ for $\theta\in (0,1)$. Previous work \cite{rosenzweig_mean-field_2022-1,rosenzweigMeanfieldApproximationHigherdimensional2022} by the author, limited to the Coulomb case, had obtained a defect corresponding to $\theta=1$. This improvement is made possible by the celebrated \emph{Brezis-Wainger-Hansson inequality} \cite{brezis_note_1980,hansson_imbedding_1979}: if $1<p<\infty$, $1\le q\le \infty$, and $aq>\d$, then 
\begin{align}\label{eq:BW}
     \forall f\in\Sc(\R^\d), \qquad \|f\|_{L^\infty} \le C\|f\|_{\dot{W}^{\frac{\d}{p},p}}\Big(1+\log^{\frac{p-1}{p}}(1+\frac{\|f\|_{\dot{W}^{a,q}}}{\|f\|_{\dot{W}^{\frac{\d}{p},p}}})\Big),
\end{align}
where $C=C(\d,p,a,q)>0$. Remark that the exponent $\frac{p-1}{p}$ is sharp. 
Brezis and Wainger \cite[Theorem 1]{brezis_note_1980} first proved the estimate \eqref{eq:BW}, which is a  generalization of an earlier $\d=p=2$ inequality due to Brezis and Gallou\"et \cite{brezis_nonlinear_1980}. 
Independently, Hansson \cite{hansson_imbedding_1979}  proved an equivalent inequality expressed in potential-theoretic language, and an alternative proof of \eqref{eq:BW} is also given in \cite{engler_alternative_1989}. 
The  inequality \eqref{eq:BW} is a substitute for the false endpoint Sobolev embedding $\|f\|_{L^\infty} \lesssim \|f\|_{\dot{W}^{\frac{\d}{p},p}}$, except when $\d=p=1$. 
In some sense, it is also dual to  the famous Trudinger-Moser inequality \cite{trudinger_imbeddings_1967, moser_sharp_1970} (see \cite[Theorem A]{brezis_note_1980}).


\section{Mean-field limits}\label{sec:MF}
We transition to applications of \cref{thm:FIprime}, beginning with sharp rates of convergence for the mean-field dynamics of Riesz gases. This section draws on the works \cite{rosenzweig_sharp_nodate,hess-childs_sharp_2025, hess-childs_sharp_2025} joint with Hess-Childs and Serfaty. 

\subsection{Background and motivation}\label{ssec:MFbackground}
In recent years, the study of {mean-field limits} for interacting particle systems has seen significant progress. A basic model is the first-order system%
	\begin{equation}\label{eq:MFode}
			d{x}_i^t = \displaystyle\frac{1}{N}\sum_{1\leq j\leq N : j\neq i}\M\nabla_x\g(x_i^t,x_j^t)dt - \mathsf{V}^t(x_i^t)dt + \sqrt{2/\be}dW_i^t, \qquad  i\in [N],
	\end{equation}
where $\be \in (0,\infty]$ has the interpretation of inverse temperature, $W_1,\ldots,W_N$ are iid $\d$-dimensional Wiener processes, and the stochastic differential $dW_i$ is understood in the It\^{o} sense. Here, $\g$ is an \emph{interaction potential} to be specified, and $\mathsf{V}^t$ is a (possibly time-dependent) \emph{external field}. $\M$ is a $\d\times\d$ matrix with real entries that specifies the type of dynamics, with $\M=-\I$ (\emph{gradient}/dissipative) and $\M$ antisymmetric (\emph{Hamiltonian}/conservative) being the most common. We assume that
	\begin{equation}\label{hofis:eq:Mnd}
	\M\xi\cdot\xi \leq 0 \qquad \forall \xi\in\R^\d,
	\end{equation}
 which is a repulsivity assumption. One can check that our assumption \eqref{hofis:eq:Mnd} for $\M$ ensures that the energy for \eqref{eq:MFode} is almost surely (a.s.) nonincreasing, therefore if the particles are initially separated, they a.s. remain separated for all time, so that there is a unique, global strong solution to the system \eqref{eq:MFode}.
 
 The mean-field limit refers to the convergence as $N \to \infty$ of the {empirical measure} $\mu_N^t$
associated to a solution $\ux_N^t$ of the system \eqref{eq:MFode}. Assuming the initial points $\ux_N^0$ are such that $\mu_N^0$ converges to a sufficiently regular measure $\mu^0$, a formal calculation leads one to expect that for $t>0$, $\mu_N^t$ converges to the solution of the \emph{mean-field equation}
	\begin{equation}\label{eq:MFlim}
			\partial_t \mu= \div ((\mathsf{V}-\M \nabla \g*\mu) \mu) + \frac1\be \Delta\mu.
	\end{equation}
	Convergence of the empirical measure is qualitatively equivalent to {\it propagation of molecular chaos} (see \cite{golse_dynamics_2016,hauray_kacs_2014} and references therein). This latter notion means that if $f_N^0(x_1, \dots, x_N) = \mathrm{Law}(X_N^0)$ is $\mu^0$-chaotic (i.e.~the $k$-point marginals $f_{N;k}^0\rightharpoonup (\mu^0)^{\otimes k}$ as $N\rightarrow\infty$ for every fixed $k$), then $f_N^t = \mathrm{Law}(\ux_N^t)$ is $\mu^t$-chaotic. There is also the stronger notion of {\it generation of chaos}, discussed in \cite{rosenzweig_modulated_2025},  which expresses that the joint law becomes $\mu^t$-chaotic as $t\rightarrow\infty$. 
	
    Mean-field limits for systems of the form \eqref{eq:MFode} have a long history, beginning with regular drifts (typically, Lipschitz) \cite{Dobrushin1979, Sznitman1991} and then progressing to increasingly more singular interactions over the years. A proper review of this history is beyond the scope of these pages. The interested reader may consult the survey \cite{CD2021}, the lecture notes \cite{Golse2022ln,golse_dynamics_2016, JW2017_survey, Jabin2014, serfaty_lectures_2024}, and the introductions of \cite{serfaty_mean_2020,nguyen_mean-field_2022}. 
    Instead, we focus on the much more recent and challenging treatment of singular interactions, specifically {log/Riesz} potentials, which we consider here under the same assumptions as above, focusing on the zero temperature case $\be=\infty$. At the end, in \cref{ssec:MFposT}, we comment on the positive temperature case.

    Mean-field convergence for the full range $-2<\s<\d$ has only been resolved in the last few years: the nonsingular case  $-2<\s<0$ \cite{rosenzweig_wasserstein_nodate}, the sub-Coulomb case $\s<\d-2$, \cite{hauray_wasserstein_2009,carrillo_derivation_2014}, the Coulomb/super-Coulomb case $\d-2\leq \s<\d$ \cite{duerinckx_mean-field_2016,carrillo_mass-transportation_2012,berman_propagation_2019,serfaty_mean_2020}, and the full singular case $0\leq \s<\d$ \cite{nguyen_mean-field_2022}. 
    These advances are due in large part to the {modulated-energy} method, which we review in the next subsection. 

	{Other techniques, such as based on weighted estimates for hierarchies of cumulants or marginals, can show propagation of chaos even for some Riesz interactions \cite{BDJ2024}. This approach is better suited to the positive temperature/diffusive case, in contrast to our zero-temperature setting, where stronger results are known \cite{BJS2022, BDJ2024,duerinckxCorrelationEstimatesBrownian2025}. Although unable to treat the full {potential} Riesz range, these approaches do have advantages in terms of working for both first- and second-order dynamics, are robust to the precise form of the interaction, and in some cases achieve sharp rates \cite{Lacker2023, LlF2023, hCR2023, Wang2024sharp}, though for essentially at most log singularities. These results imply probabilistic statements about the mean-field convergence rates, but they do not imply the deterministic convergence rates obtainable with our methods.}

\subsection{Main results}\label{ssec:MFmain}
The optimal rate of convergence as $N\rightarrow\infty$ of the empirical measure $\mu_N^t$ to the solution $\mu^t$ of \eqref{eq:MFlim} in the distance $\Fr_N$ is $N^{\frac{\s}{\d}-1}$. The reason is that {minimizers} $X_N$ of Riesz interaction energies with a suitable confinement $V$ in $\R^\d$ or on the flat torus are known to satisfy $|\Fr_N(\XN,\mu_V)| \propto N^{\frac{\s}{\d}-1}$, where $\mu_V$ is the minimizer of the mean-field {energy}, called the equilibrium measure (see \cref{ssec:Lakeback}). This has been shown for the (super-)Coulomb case \cite{sandier_1d_2015,sandier_2d_2015,rougerie_higher-dimensional_2016,petrache_next_2017}, the general singular case $0\le \s<\d$ on the flat torus \cite{hardinNextOrderEnergy2017}, and for the nonsingular case $-2<\s<0$ \cite{hess-childs_optimal_nodate}. Since minimizers are critical points of the energy, it follows that when $\mathsf{V}^t = \nab V$, they yield stationary solutions of \eqref{eq:MFode} with $\M=-\mathbb{I}$ and $\be=\infty$. Similarly, ${\mu_V}$ is a stationary solution of \eqref{eq:MFlim}.  Sharp rates for the nonstationary case were only completely established recently: the Coulomb case \cite{serfaty_gaussian_2023, rosenzweig_mean-field_2022-1},  (super-)Coulomb case \cite{rosenzweig_sharp_nodate}, full singular case \cite{hess-childs_sharp_2025}, and the nonsingular case \cite{rosenzweig_wasserstein_nodate}. 

\cref{thm:mainMF} below summarizes the state of the art for mean-field convergence (at zero temperature) in the full potential Riesz case $-2<\s<\d$ with the sharp rate. The $\d-2\le \s<\d$ portion is taken from {\cite{rosenzweig_sharp_nodate}}, the $0\le \s<\d-2$ portion from {\cite{hess-childs_sharp_2025}}, and the $-2<\s<{0}$ portion from \cite{rosenzweig_wasserstein_nodate}. Remark that there is no additive $N^{\frac{\s}{\d}-1}$ error in the nonsingular case $-2<\s\le 0$ because for any suitable $\mu$, it has been shown that $\min_{\XN}\Fr_N(\XN,\mu)\propto N^{\frac{\s}{\d}-1}$, in particular nonnegative. 
As explained in the next subsection, the proof is a quick consequence of the by now standard modulated-energy method and \cref{thm:FIprime}. Following the statement, we record some remarks on the theorem's assumptions and implications.

\begin{thm}\label{thm:mainMF}
Let $\g$ be of the form \eqref{eq:gmod}, and suppose that for some $\as\in (\d,\d+2)$ and $T>0$, $\mathsf{V}$ satisfies
\begin{align}
\int_0^T\big(\|\nabla\mathsf{V}^t\|_{L^\infty} + \|\Dm^{\frac{\as}{2}}\nab \mathsf{V}^t\|_{L^{\frac{2\d}{\as-2}}}\indic_{\substack{\as>2 \\ 0\le \s<\d-2}} + \|\Dm^{\frac{\d-\s}{2}}\mathsf{V}^t\|_{L^{\frac{2\d}{\d-\s}}}\indic_{-2<\s<0}\big)dt<\infty.
\end{align}
Assume the equation \eqref{eq:MFlim} admits a solution $\mu \in L^\infty([0,T],\P(\R^\d)\cap L^\infty(\R^\d))$ such that $\int_{(\R^\d)^2} |\g|(x-y)d(\mu^t)^{\otimes 2}<\infty$ for every $t\in [0,T]$ and 
\begin{multline}\label{nmut}
\mathcal{N}(u^T) \coloneqq \int_0^T \big(\|\nab u^t\|_{L^\infty} + \|\Dm^{\frac{\as}{2}}u^t\|_{L^{\frac{2\d}{\as-2}}}\indic_{\substack{\as>2 \\ 0\le \s<\d-2}} + \|\Dm^{\frac{\d-\s}{2}}u^t\|_{L^{\frac{2\d}{\d-\s}}}\indic_{-2<\s<0}\big)dt< \infty, \\
 \text{where} \ u^t\coloneqq -\M\nab\g\ast\mu^t +\mathsf{V}^t.
\end{multline}
There exists a constant  $C>0$, depending only on $\d$, $\s$, such that for any solution $\ux_N$ of \eqref{eq:MFode},
\begin{multline}\label{distcoul}
\Fr_N(\ux_N^t,\mu^t) + \frac{\log (N\|\mu^t\|_{L^\infty}) }{2 N\d} \indic_{\s=0} + C\|\mu^t\|_{L^\infty}^{\frac{\s}{\d}}N^{\frac{\s}{\d}-1}\\
 \le Ce^{C\mathcal{N}(u^t)}\Bigg(\Fr_N(\ux_N^0,\mu^0) + \sup_{\tau \in [0,t]}\Big( \frac{\log (N\|\mu^\tau\|_{L^\infty}) }{2 N\d} \indic_{\s=0} + C\|\mu^\tau\|_{L^\infty}^{\frac{\s}{\d}}N^{\frac{\s}{\d}-1}\Big)\Bigg).
\end{multline}
Consequently, if $\lim_{N\to \infty}F_N(\ux_N^0, \mu^0) =0$, then $\mu_N^t \rightharpoonup \mu^t$ for every $t\in [0,T]$.
\end{thm}


{
\begin{remark}\label{rem:V=0}
When $\mathsf{V}^t=0$, the norm $\|\mu^t\|_{L^\infty}$. 
Consequently,
\begin{align}
\sup_{\tau \in [0,t]}\Big( \frac{\log (N\|\mu^\tau\|_{L^\infty}) }{2 N\d} \indic_{\s=0} + C\|\mu^\tau\|_{L^\infty}^{\frac{\s}{\d}}N^{\frac{\s}{\d}-1}\Big) \le \frac{\log (N\|\mu^0\|_{L^\infty}) }{2 N\d} \indic_{\s=0} + C\|\mu^0\|_{L^\infty}^{\frac{\s}{\d}}N^{\frac{\s}{\d}-1}.
\end{align}
\end{remark}
}

{
\begin{remark}\label{rem:extVFreg}
The reader may check that the above condition for $\mathsf{V}^t = -c^t x + \Phi^t(x)$ is satisfied, where $\int_0^T |c^t|dt <\infty$, and  $\int_0^T(\|\nab \Phi^t\|_{L^\infty} + \|\Dm^{\frac{\as}{2}} \Phi^t\|_{L^{\frac{2\d}{\as-2}}}\indic_{\as>2})dt < \infty$. 
The regularity assumption \eqref{nmut} for $\mu^t$ follows from a number of works of in the PDE literature. We refer to \cite[Section 1.3]{serfaty_mean_2020} for references in the (super-)Coulomb case and \cite[Section 5.2]{hess-childs_sharp_2025} for direct verification in the sub-Coulomb case. 
\end{remark}
}

{
\begin{remark}\label{rem:pc}
By a standard argument (see, e.g., \cite[Remark 3.7]{serfaty_mean_2020} or \cite[Remark 1.5]{rosenzweig_global--time_2023} for details), \cref{thm:mainMF} implies propagation of chaos for the system \eqref{eq:MFode}. However, this ``global-to-local'' argument in general leads to a suboptimal rate  \cite{Lacker2023}.
\end{remark}
}

\subsection{Modulated energy and commutators}\label{ssec:MFcomm}
    Since the empirical measure $\mu_N^t$ associated with the microscopic system ~\eqref{eq:MFode} is a weak solution of the mean-field equation ~\eqref{eq:MFlim}, a natural approach for showing mean-field convergence is a weak-strong stability estimate for~\eqref{eq:MFlim}. For Riesz interactions, such stability is conveniently formulated in terms of the modulated energy due to its coercivity---vanishing of the modulated energy metrizes the convergence of $\mu_N^t$ to $\mu^t$.

The modulated energy approach to mean-field convergence consists of establishing a Gr\"onwall relation for $\Fr_N(\XN^t,\mu^t)$, where $\XN^t$ is a solution of the microscopic system \eqref{eq:MFode} and $\mu^t$ is a solution of the mean-field equation \eqref{eq:MFlim}.  The key observation  is that  $\Fr_N(\XN^t,\mu^t)$ satisfies the dissipation inequality (see, e.g., \cite[Section 3.1]{serfaty_mean_2020})
\begin{align}\label{eq:modenergyineq}
\frac{d}{dt}\Fr_N(\XN^t,\mu^t) \le \int_{(\R^\d)^2\setminus \triangle} \nab\g(x-y)\cdot (u^t(x)-u^t(y)) d(\mu_N^t-\mu^t)^{\otimes 2},
\end{align}
where $u^t$ is the mean-field velocity from \eqref{nmut}. The right-hand side is exactly a first-order commutator expression with $v=u^t$. Importantly, the transport only depends on the mean-field density $\mu^t$, not the empirical measure $\mu_N^t$. The beauty of the modulated-energy method is that with this inequality in hand, one can apply the estimates of \cref{thm:FIprime} and conclude by Gr\"onwall's lemma. For the exact details on closing the differential inequality and obtaining the estimate \eqref{distcoul} stated in \cref{thm:mainMF}, we refer to \cite[Section 6]{rosenzweig_sharp_nodate} or \cite[Section 5.1]{hess-childs_sharp_2025}.



\subsection{Optimal regularity assumptions}\label{ssec:MFreg}
The assumption that $\mu^t\in L^\infty$ is not so important and can be greatly relaxed at the cost of worsening the rate of convergence. More fundamental is the condition \eqref{nmut} in the statement of \cref{thm:mainMF}, needed to apply \cref{thm:FIprime}, which enforces some regularity for $\mu^t$. Given the relationship between the scaling-critical regularity and that for which uniqueness is expected to hold for \eqref{eq:MFlim}, it is natural to ask if the result holds when $\mu^t$ belongs to the critical space $\dot{W}^{\frac{\d}{p}+\s+2-\d,p}$ for $1<p<\infty$. Provided that $p \le \frac{2\d}{\d-\s-2}$, one has $\|\Dm^{\frac{\d-\s}{2}}u^t\|_{L^{\frac{2\d}{\d-\s-2}}} \lesssim \|\mu^t\|_{\dot{W}^{\frac{\d}{p}+\s+2-\d,p}}$. {But in general}, $\nab u^t$ is only in BMO when $\mu^t\in \dot{W}^{\frac{\d}{p}+\s+2-\d,p}$. As we observed in \cite{rosenzweig_mean-field_2022-1,rosenzweigMeanfieldApproximationHigherdimensional2022}, even in this borderline regularity case, one can still obtain a Gr\"onwall-type\footnote{More precisely, an estimate in the form of Osgood's lemma (see, e.g., \cite[Section 3.1.1]{bahouri_fourier_2011}).}  estimate for the modulated energy, which is weaker than \eqref{distcoul} but still suffices for proving mean-field convergence. Tying back into the discussion of \cref{ssec:rcomtransreg}, the key ingredient was a defective commutator estimate. These {works} \cite{rosenzweig_mean-field_2022-1,rosenzweigMeanfieldApproximationHigherdimensional2022} considered only the Coulomb case, where the natural scale-invariant space is $L^\infty$, and had a short-time restriction for $\d\ge 3$.

In \cite[Theorem 2.8]{hess-childs_sharp_2025}, we used our new defective commutator estimate to generalize the results of \cite{rosenzweig_mean-field_2022-1,rosenzweigMeanfieldApproximationHigherdimensional2022} to the entire Riesz case without any restrictions on the interval of time for which mean-field convergence holds. In the sub-Coulomb case, we showed that if the initial mean-field density $\mu^\circ \in L^{\frac{\d}{\d-\s-2}}$, which is scaling-critical, then {there is} a unique global solution to the mean-field equation \eqref{eq:MFlim} which is a (quantitative) mean-field limit of the microscopic system \eqref{eq:MFode}. For details on {closing} the resulting differential inequality for the modulated energy after applying the defective commutator estimate, we refer to \cite[Sections 2.3, 4.4]{hess-childs_sharp_2025}.

\subsection{Positive temperature}\label{ssec:MFposT}
So far, we have only discussed the zero-temperature case $\be=\infty$. At positive temperature $\be<\infty$, where the dynamics of \eqref{eq:MFode} are now stochastic, the pure modulated-energy approach described above can be extended, considering now moments of the modulated energy, at least in the {sub-Coulomb} case $\s<\d-2$. This restriction is due to the It\^o correction to the evolution of the energy, which is not clear how to suitably control when $\s\ge \d-2$ (barely missing the Coulomb case). We refer to \cite{rosenzweig_global--time_2023,hess-childs_large_2023,hess-childs_sharp_2025} for details on this extension and the resulting convergence rates.

{ For gradient dynamics (i.e. $\M=-\I$) at positive temperature, Bresch et al. \cite{bresch_mean-field_2019, bresch_modulated_2019, bresch_mean_2023} beautifully observed that the modulated energy combines naturally with the previously used relative entropy \cite{jabin_quantitative_2018, guillin_uniform_2024, feng_quantitative_2024, rosenzweig_relative_2024} in the form of the \emph{modulated free energy}
\begin{align}\label{eq:mfe}
    \mathsf{E}_N(f_N^t,\mu^t)&\coloneqq \frac{1}{\be N}\mathsf{H}\left(f_N^t\mid (\mu^t)^{\otimes N}\right)+\E_{X_N\sim f_N^t}\Big[\Fr_N(\XN,\mu^t)\Big]     = \frac{1}{N}\mathsf{H}(f_N^t \mid \Q_{N,\be}(\mu^t)) - \frac{\log \mathsf{K}_{N,\be}(\mu^t)}{N},
\end{align}
where $\mathsf{H}$ is the relative entropy/KL divergence. The second equality is the reformulation in terms of the modulated Gibbs measure and partition function introduced by the author and Serfaty \cite{rosenzweig_modulated_2025},
\begin{align}\label{eq:modGibbs}
    d\mathbb{Q}_{N,\be}(\mu)\coloneqq \frac{e^{-\be N\Fr_N(\XN,\mu)}}{\mathsf{K}_{N,\be}(\mu)}\,d\mu^{\otimes N}, \qquad \mathsf{K}_{N,\be}(\mu)  \coloneqq \int_{(\R^\d)^N}e^{-\be N\Fr_N(\XN,\mu)}\,d\mu^{\otimes N}(\XN).
\end{align}
The quantity \eqref{eq:mfe} is well-suited to proving entropic propagation of chaos \cite{bresch_mean-field_2019, bresch_modulated_2019, bresch_mean_2023, chodron_de_courcel_sharp_2023, rosenzweig_modulated_2025, rosenzweig_relative_2024, cai_propagation_2024} and can even handle logarithmically attractive interactions \cite{bresch_mean-field_2019, bresch_mean_2023, chodron_de_courcel_attractive_2025}.  Commutator estimates play the same essential role in the positive temperature setting.
More precisely, if $f_N^t\coloneqq \mathrm{Law}(\XN^t)$, the original BJW quantity can be written in our notation as
}

\section{Supercritical mean-field limits}\label{sec:Lake}
We conclude with the second application of our sharp commutator estimates to  so-called supercritical mean-field limits. Namely, these estimates combined with a modulated-energy scheme allow to derive the Lake equation from the Newtonian $N$-body problem in the joint mean-field and quasineutral limit for (super-)Coulomb Riesz interactions under optimal scaling assumptions. The material in this section is drawn from the joint work \cite{rosenzweig_lake_2025} with Serfaty.  

\subsection{Background and motivation}\label{ssec:Lakeback}
Consider a \emph{Newtonian system} of $N$ particles with a pairwise {interaction potential} $\g$ and external {confining potential} $V$:
\begin{equation}\label{lake:eq:NewODE}
\begin{cases}
\dot{x}_i^t = v_i^t \\
\dot{v}_i^t = -\ga v_i^t \displaystyle-\frac{1}{\vep^2 N} \sum_{1\leq j\leq N: j\neq i}\nabla\g(x_i^t-x_j^t) -\frac{1}{\vep^2}\nabla V(x_i^t)\\
(x_i^0,v_i^0) = (x_i^\circ,v_i^\circ),
\end{cases}
\qquad 1\le i\le N.
\end{equation}
The positions and velocities are assumed to belong to $\R^\d$. Here, $\ga\geq 0$ is the \emph{friction coefficient}, and $\vep>0$ is a small parameter, possibly depending on $N$, which encodes physical information about the system.  In contrast to the first-order dynamics  in \eqref{eq:MFode}, the dynamics in \eqref{lake:eq:NewODE} are now second-order. For concreteness and because are results will ultimately be for this class, we will limit our attention to (Coulomb/super-Coulomb) Riesz potentials. Though, much of the discussion at this stage makes sense for more general interaction potentials.

{We are interested in the large $N$ and small $\vep$ limit of \eqref{lake:eq:NewODE}, the latter of which we interpret as a \emph{quasineutral limit} elaborated on in the next subsection. There are two (mathematically equivalent) motivations for our setup.

The first motivation is non-neutral plasmas \cite{OD1998} (see also \cite{WBIP1985,MKTL2008} for relevance to trapped neutral systems). The system \eqref{lake:eq:NewODE} models the evolution of a trapped system of ions near thermodynamic equilibrium, meaning the spatial density $\mu_N^t \coloneqq \frac1N\sum_{i=1}^N\delta_{x_i^t}$ is close to the \emph{equilibrium measure} $\mu_V$. This equilibrium measure is defined as
the probability measure that minimizes the {macroscopic}  or mean-field energy
\begin{equation}\label{lake:eq:Edef}
\mathcal{E}(\mu) \coloneqq \int_{\R^\d}Vd\mu + \frac{1}{2}\int_{(\R^\d)^{ 2}}\g(x-y)d\mu(x) d\mu(y).
\end{equation}
We refer to \cite[Chap. 2]{serfaty_lectures_2024} for  details on the equilibrium measure in this Coulomb/Riesz context, as well as for its connection to the solution of the fractional obstacle problem, which will be also discussed below. 
One can see the equilibrium measure as a generalization of  the uniform measure on a torus when considering an infinitely extended trapped system.

The second motivation comes from two-species globally neutral systems. In this setting, the empirical spatial density $\mu_N^t$ is close to a fixed density $\mu$, representing the density of a stationary background of an oppositely charged species of particles (e.g.~heavy positively charged ions). Alternatively, one may think of this as a one-component plasma with a nonuniform background. In this case, the term $-\nab V(x_i^t)$ on the right-hand side of the second line of \eqref{lake:eq:NewODE} should be replaced by the attractive force $+\nab(\g\ast\mu)(x_i^t)$ due to the background.

The settings for each of these motivations are mathematically equivalent. The latter corresponds to $V = -\g\ast\mu + \frac12\int_{(\R^\d)^2}\g(x-y)d\mu^{\otimes 2}$. While the former corresponds to considering a background density such that $\nabla(\g\ast\mu+V)=0$ in the support of $\mu$. This is in particular achieved when (but not only when) $\mu$ is equal to the equilibrium measure $\mu_V$. 
}


Under suitable assumptions on the external potential $V$, our goal is  to show that if the initial \emph{empirical spatial density} $\mu_N^\circ \coloneqq \frac1N\sum_{i=1}^N \delta_{x_i^\circ} \rightharpoonup \mu_V$ as $N\rightarrow\infty$ and the initial velocities $v_i^\circ \approx u^\circ(x_i^\circ)$, for a macroscopic vector field $u^\circ$ on $\R^\d$, then the total {empirical measure} $\frac{1}{N}\sum_{i=1}^N \delta_{(x_i^t,v_i^t)}$ associated to a solution of \eqref{lake:eq:NewODE} converges as $\vep\rightarrow 0$ and $N\rightarrow\infty$ to the \emph{monokinetic} measure $\mu_V(x)\delta_{u^t(x)}(v)$, where $u^t$ satisfies the \emph{Lake equation}
\begin{equation}\label{lake:eq:Lake}
\begin{cases}
\p_t u +\ga u+ u \cdot\nabla u = -\nabla p\\
\div(\mu_V u) = 0.
\end{cases}
\end{equation}
In particular, when $\g$ is the Coulomb potential, we give a microscopic counterpart to the proof of the \emph{quasineutral} limit for Vlasov-Poisson with monokinetic data by  Barr\'{e} et al. \cite{BCGM2015}.

Note that if $\mu_V$ is constant and $\ga=0$, then \eqref{lake:eq:Lake} is nothing but the \emph{incompressible Euler equation}. The \emph{pressure} $p$ is a Lagrange multiplier to enforce the incompressibility constraint $\div(\mu_V u)=0$. Multiplying both sides of the first equation of \eqref{lake:eq:Lake} by $\mu_V$ and taking the divergence, the pressure $p$ is obtained from the velocity $u$ by solving the divergence-form elliptic equation
\begin{equation}\label{lake:eq:press}
-\div\pa*{\mu_V\nabla p} = \div^{2}\pa*{\mu_V u^{\otimes 2}}= \div\pa*{\mu_V u\cdot\nabla u}.
\end{equation}

{Equation \eqref{lake:eq:Lake}, which is also sometimes called the \emph{anelastic equation}, appears in the modeling of atmospheric flows \cite{OP1962, Masmoudi2007} and superconductivity \cite{CR1997, DS2018, Duerinckx2018} and has been mathematically studied in \cite{LOT1996, LOT1996phy, BCGM2015, Duerinckx2018}.  The equation also arises as a mean-field limit for Ginzburg-Landau vortices with pinning and forcing  \cite{DS2018} have shown that the equation. {Additionally, the Lake equation has been shown \cite{Menard2023lake} to be a mean-field limit for a model of vortices in shallow water with varying topography \cite{Richardson2000}.}}



\subsection{The combined mean-field and quasineutral limit}\label{lake:ssec:introMFQN}
To see how the equation \eqref{lake:eq:Lake} appears as a formal limiting dynamics for the empirical measure of \eqref{lake:eq:NewODE}, we argue as follows.

Suppose that the parameter $\vep>0$ is fixed. Analogous to the derivation of the equation \eqref{eq:MFlim} in \cref{ssec:MFbackground}, the mean-field limit in this context refers to the convergence as $N\rightarrow\infty$ of the total empirical measure $f_{N,\vep}^t \coloneqq \frac{1}{N}\sum_{i=1}^N \delta_{z_i^t}$, where $z_i^t := (x_i^t,v_i^t)$,   to a solution $f_\vep^t$ of the \emph{Vlasov equation with friction} 
\begin{equation}\label{lake:eq:Vlas}
\begin{cases}
\p_t f_\vep+v\cdot\nabla_x f_\vep-\frac{1}{\vep^2}\nabla(V+\g\ast\mu_\vep)\cdot\nabla_v f_\vep - \div_v(\ga vf)=0\\
\mu_\vep = \int_{\R^\d}df_\vep(\cdot,v) \\
f_\vep|_{t=0} = f_\vep^\circ.
\end{cases}
\end{equation}
When $\g$ is Coulomb and $\ga=0$, the equation \eqref{lake:eq:Vlas} is known as \emph{Vlasov-Poisson}. More generally, for $\g$ as in \eqref{eq:gmod}, the equation is called \emph{Vlasov-Riesz}. It is a difficult problem to derive the Vlasov-Poisson/Vlasov-Riesz equation directly from \eqref{lake:eq:NewODE}. While the case of regular potentials $\g$ (e.g. globally $C^{1,1}$) \cite{NW1974, BH1977, Dobrushin1979, Duerinckx2021gl} or even just potentials with bounded force $\nabla\g$ \cite{JW2016} is well understood, the Coulomb case in general remains out of reach except in dimension 1 \cite{Trocheris1986,Hauray2014}. To our knowledge, the best results for singular potentials are limited to forces $\nabla\g$ which are square integrable at the origin \cite{HJ2007, HJ2015, BDJ2024} (which barely misses the 2D Coulomb case) or are for Coulomb potentials with short-distance vanishing cutoff \cite{BP2016, Lazarovici2016, LP2017, Grass2021,feistl-heldMeanfieldLimitVlasovPoisson2025,feistl-heldMeanfieldLimitVlasovPoisson2025a}.\footnote{If one adds noise to the velocity equation in \eqref{lake:eq:NewODE}, corresponding to the \emph{Vlasov-Fokker-Planck} mean-field limit, then the 2D Coulomb case has been shown \cite{BJS2022}. See also \cite{duerinckxCorrelationEstimatesBrownian2025} for some Riesz cases.} Recently, Duerinckx and Serfaty \cite[Appendix]{serfaty_mean_2020} proved the mean-field limit for the Vlasov-Poisson equation---and more generally super-Coulomb Vlasov-Riesz\footnote{Combining the total modulated energy introduced in this work with the commutator estimate \cref{thm:FIprime}\eqref{eq:FIsubC1}, one can extend this result to the sub-Coulomb Vlasov-Riesz equation as well.}---for monokinetic/cold initial data, for which the Vlasov-Poisson equation reduces to the \emph{pressureless Euler-Poisson equation}. Putting aside this question of rigorously showing the mean-field limit, we seek to {formally}  derive the Lake equation \eqref{lake:eq:Lake} from the Vlasov equation \eqref{lake:eq:Vlas} in the limit as $\vep\rightarrow0$.

The regime under consideration is when the spatial density $\mu_\vep^t$ converges to the equilibrium measure $\mu_V$ as $\vep\rightarrow 0$ (this is an assumption). Decomposing the potential 
\begin{equation}
V+\g\ast\mu_\vep = (V+\g\ast\mu_V) + \g\ast(\mu_\vep-\mu_V),
\end{equation}
the fact that, {by characterization of the equilibrium measure}, $V+\g\ast\mu_V$ is constant on the support of $\mu_V$ implies that
\begin{equation}
\nabla\pa*{V+\g\ast\mu_\vep} = \nabla\g\ast(\mu_\vep-\mu_V), \qquad x\in\supp \mu_V.
\end{equation}
Assuming that the renormalized electric potential difference $\frac{1}{\vep^2}\g\ast(\mu_\vep-\mu_V)$ has a weak limit $p$ as $\vep\rightarrow 0$, we see that the weak limit $f\coloneqq \lim_{\vep\rightarrow 0}f_\vep$ satisfies the equation
\begin{equation}\label{lake:eq:KIE}
\begin{cases}
\p_tf+v\cdot\nabla_xf - \nabla p\cdot\nabla_vf -\div_v (\ga v f)=0\\
\mu_V = \int_{\R^\d}df(\cdot,v)\\
f|_{t=0} = f^\circ,
\end{cases}
\end{equation}
where $f^\circ$ is the weak limit of $f_\vep^\circ$.

Let us now define the \emph{current} $J(x)\coloneqq \int_{\R^\d}v\, df(x,v)$ associated to \eqref{lake:eq:KIE}. Integrating both sides of the first equation in \eqref{lake:eq:KIE} with respect to {$v$}, then using that the spatial density is equal to $\mu_V$ for all time, we find that {$\div J^t$ is constant in time}. After some calculus (see \cite[Section 1.1]{rosenzweig_lake_2025} for the computation), one finds that
\begin{equation}\label{lake:eq:Jeq}
\p_t J + \div\int_{\R^\d}v^{\otimes 2}df(\cdot,v) = -\ga J -\mu_V\nabla p.
\end{equation}
This equation is not closed in terms of $(J,p)$, since the second term on the left-hand side requires knowledge of the second velocity moment of $f$, {which in turn depends on third moment and so on} (this is the famous closure problem for moments of the Vlasov equation; see, e.g., \cite{Uhlemann2018}). But making the \emph{monokinetic} or \emph{``cold electrons''}\footnote{This terminology from in the physics literature stems from the fact that the temperature of the distribution is zero.} ansatz $f(x,v)=\mu_V(x)\delta(v-u(x))$, it follows that $J=\mu_V u$. Since $\mu_V$ is independent of time, substituting this identity into \eqref{lake:eq:Jeq} yields
\begin{equation}\label{lake:eq:LakeJ}
\begin{cases}
\mu_V\p_t u +\div(\mu_V u^{\otimes 2}) = -\mu_V\pa*{\ga u+\nabla p},\\
\div(\mu_V u) = 0.
\end{cases}
\end{equation}
Assuming that $\mu_V$ is positive on its support, we see that \eqref{lake:eq:LakeJ} is equivalent to  \eqref{lake:eq:Lake}.

{In the plasma physics setting of Vlasov-Poisson, the limit $\vep\rightarrow 0$ is called the \emph{quasineutral limit}, and the equation \eqref{lake:eq:KIE}, in the case $\mu_V\equiv 1$, is called the \emph{kinetic incompressible Euler (KIE) equation} \cite{Brenier1989}. The inhomogeneous case of \eqref{lake:eq:KIE} does not seem to have previously appeared in the literature. Nowhere in the above reasoning did we assume a specific form for $\g$ (e.g. Coulomb). This demonstrates a certain \emph{universality} of the KIE for this kind of singular limit, which appears to be new. In this plasma physics setting, the distribution function $f$ models the evolution of electrons against a stationary background of positively charged ions. After a rescaling to dimensionless variables, the parameter $\vep$ corresponds to the \emph{Debye (screening) length} of the system, which is the scale at which charge separation in the plasma occurs. When the Debye length is much smaller than the length scale of observation, the plasma is said to be quasineutral, since it appears neutral to an observer. The rigorous justification of the quasineutral limit is a difficult problem  and has been studied in \cite{BG1994, Grenier1995, Grenier1996,Grenier1999,Brenier2000,Masmoudi2001, HkH2015, HkR2016, HkI2017, HkI2017jde, GpI2018, GpI2020sing}. We refer  to the survey \cite{GpI2020qn} and references therein for further discussion.}

The preceding formal calculations suggest that in the limit as $N +\vep^{-1}\rightarrow\infty$, which one can physically interpret as a combined mean-field and quasineutral limit, the empirical measure $f_{N,\vep}^t$ of the Newtonian system \eqref{lake:eq:NewODE} converges to a solution $f^t$ of the KIE \eqref{lake:eq:KIE}, which reduces to the Lake equation \eqref{lake:eq:Lake} for monokinetic solutions. Thus, we expect that if the particle velocities $v_i^t \approx u^t(x_i)$, then the empirical measure $f_{N,\vep}^t$ converges to the measure $\delta_{u^t(x)}(v)\mu_V(x)$ as $N+\vep^{-1}\rightarrow\infty$, where $u^t$ is a solution of \eqref{lake:eq:LakeJ}. A more general interpretation of the limit as $N +\vep^{-1}\rightarrow\infty$ is as a \emph{supercritical mean-field limit} of the system \eqref{lake:eq:NewODE}. This terminology coined by Han-Kwan and Iacobelli \cite{han-kwan_newtons_2021} refers to the fact that the force experienced by a single particle in \eqref{lake:eq:NewODE} formally diverges as $N\rightarrow\infty$, compared to being $O(1)$ for the usual $1/N$ mean-field scaling. 

\subsection{Scaling assumptions and prior work}\label{lake:ssec:introPW}

{The convergence of $f_{N,\vep}^t$ to $\delta_{u^t(x)}(v)\mu_V(x)$ was previously shown by Han-Kwan and Iacobelli \cite{han-kwan_newtons_2021} in the spatially periodic Coulomb case (i.e.~$x\in\T^\d$) when $V=0$ and $\mu_V\equiv 1$ assuming $\vep N^{\frac{1}{\d(\d+1)}} \rightarrow \infty$ as $\vep+N^{-1}\rightarrow 0$, where $u^t$ is a solution of the incompressible Euler equation. In this work, they recognized (building on the aforementioned work \cite[Appendix]{serfaty_mean_2020} for the mean-field limit) that the modulated-energy method may be also used for this supercritical limit provided one adds a suitable $O(\vep^2)$ corrector to the background spatial density in the definition of the modulated potential energy. We elaborate more on this idea in \cref{lake:ssec:introMPf} below. Their proof may be viewed as a generalization of Brenier's \cite{Brenier2000} modulated-energy approach to proving the quasineutral limit of Vlasov-Poisson with monokinetic data to allow (via renormalization) for the solution of Vlasov-Poisson to be a sum of Diracs. As in all modulated-energy approaches, the key tool is the commutator estimate \eqref{eq:introcomm} and sharper additive errors allow to treat a large range of scaling relations for $\vep,N$.} 

{In the same setting, this convergence was subsequently improved by the  author \cite{rosenzweig_rigorous_2023} to
\begin{equation}\label{lake:eq:RoseSA}
\lim_{\vep+N^{-1}\rightarrow 0} \vep^{-2}N^{-\frac{2}{\d}} =0
\end{equation}
using a sharp commutator estimate for the Coulomb case \cite{leble_fluctuations_2018, serfaty_gaussian_2023, rosenzweig_rigorous_2023}. Sharper estimates in terms of the dependence on the solution $u$ of Euler's equation were also shown. \cite{rosenzweig_rigorous_2023} conjectured that the scaling assumption \eqref{lake:eq:RoseSA} should be in general optimal, in the sense that the incompressible Euler equation should not be the limiting evolution of the empirical measure when $\vep^{-2} \ll N^{-\frac2\d}$. 
}

{These previous works are limited to the Coulomb case on the torus, which is an idealized setting since it assumes a spatially uniform density. Moreover, they left open the question of a rigorous justification of the universality of the Lake equation with respect to the interaction and confinement in the sense that it only depends on the equilibrium measure. Finally, these previous works also left open the optimality of the scaling relation between $\vep$ and $N$.}


\subsection{Informal statement of main results}\label{lake:ssec:introMR}
In  \cite{rosenzweig_lake_2025} joint with Serfaty, we settled the above questions for (super-)Coulomb Riesz interactions. To present the main result of that work, we introduce the \emph{total modulated energy}
\begin{equation}\label{lake:eq:soMEdef}
\Hr_{N}(\uz_N^t, u^t) \coloneqq \frac{1}{2N}\sum_{i=1}^N |v_i^t-u^t(x_i^t)|^2 + \frac{1}{\vep^2}\Fr_N(\ux_N^t,\mu_V+\vep^2\Uu^t)  + \frac{1}{\vep^2 N}\sum_{i=1}^N \zeta(x_i^t).
\end{equation}
Here, $\uz_N^t$ is a solution of the $N$-particle system \eqref{lake:eq:NewODE}. The vector field $u^t$ is not a solution of the Lake equation \eqref{lake:eq:Lake} but rather an extension from $\supp \mu_V $ to all of $\R^\d$, such that the regularity is preserved and $\div(\mu_V u) = 0$ (see \cite[Section 4]{rosenzweig_lake_2025}). The physical Cauchy problem for \eqref{lake:eq:Lake} is in the domain $\supp\mu_V$ subject to a no-flux boundary condition, but the  microscopic dynamics $Z_N^t$ are not confined to $\supp\mu_V$. The extension allows us to compare between the two settings. 

The first term of \eqref{lake:eq:soMEdef} is the \emph{modulated kinetic energy}. The second term is the \emph{modulated potential energy} previously seen in \eqref{hofis:eq:modenergy}. 
{Here, $\Uu^t$ is an $N$-independent \emph{corrector}, whose precise definition the reader may find in \cite[Eq. (5.2)]{rosenzweig_lake_2025} (see also Remark 5.1 there).}
The third term is
\begin{align}\label{lake:eq:introzetadef}
\zeta\coloneqq \g\ast\mu + V - c
\end{align}
for Robin constant $c$, {which is nonnegative and vanishes on $\supp \mu_V$}.  The quantity \eqref{lake:eq:soMEdef} depends on both $\vep,N$; but since we view $\vep=\vep(N)$, we omit the $\vep$ dependence  to lighten the notation.

Because the precise, rigorous version requires stating some rather cumbersome regularity assumptions, we only present an informal version of the main result in \cref{lake:thm:mainSMF} below. The interested reader may find the precise statment in \cite[Theorem 5.2]{rosenzweig_lake_2025} and also a generalization of our findings to the case of \emph{regular} interactions in Appendix A of the cited work.

\begin{thm}[Informal]\label{lake:thm:mainSMF}
Suppose that the equilibrium measure $\mu_V$ is sufficiently regular and $\Sigma\coloneqq\supp  \mu_V$ {is a sufficiently smooth domain}. Assume $u$ is a sufficiently regular solution to \eqref{lake:eq:Lake} on $[0,T]$. Then there exist continuous functions $C_1,C_2,C_3: [0,T]\rightarrow \R_+$, which depend on $\d,\s, \ga$, and norms of $u$, such that for any solution $Z_N^t = (X_N^t, V_N^t)$ of \eqref{lake:eq:NewODE}, it holds that
\begin{align}\label{lake:eq:mainSMF}
\Hr_{N}(\uz_N^t, u^t) + \frac{\log N}{2\d N\vep^2}\indic_{\s=0} \leq e^{C_1^t}\Big(\Hr_{N}(\uz_N^0, u^0)   +\frac{\log N}{2\d N\vep^2}\indic_{\s=0} + \frac{C_2^tN^{\frac{\s}{\d}-1}}{\vep^2} {+ C_3^t\vep^2} \Big).
\end{align}
In particular, if 
\begin{equation}\label{lake:eq:SMFscl}
\lim_{\vep+\frac{1}{N}\rightarrow 0} \Big(\Hr_{N}(\uz_N^0,u^0) + \frac{\log N}{2\d N\vep^2}\indic_{\s=0}\Big) =  \lim_{\vep+\frac{1}{N}\rightarrow 0}  \frac{N^{\frac{\s}{\d}-1}}{\vep^2} =0,
\end{equation}
then we have the weak-* convegence
\begin{equation}\label{lake:eq:SMFweak}
\frac1N\sum_{i=1}^N\delta_{z_i^t} \xrightharpoonup[\vep+\frac{1}{N}\rightarrow 0]{} \delta_{u^t(x)}(v)\mu_V(x), \qquad \forall t\in [0,T].
\end{equation}
\end{thm}

{The scaling assumption $N^{\frac{\s}{\d}-1}/\vep^2\rightarrow 0$ in \eqref{lake:eq:SMFscl} is in general optimal, in the sense that there exists a sequence of solutions $\uz_N^t$ to \eqref{lake:eq:NewODE} such that ${\frac1N\sum_{i=1}^N \delta_{z_i^t}} \xrightharpoonup[]{} \delta_{u^t(x)}(v)\mu_V(x)$ as $\vep+\frac1N\rightarrow 0$, but the total modulated energy $\Hr_{N}(\uz_N^t,u^t)$ does not vanish. This follows from essentially the same reasoning explained in \cref{sec:MF} for the sharp rate of convergence for the mean-field limit in the modulated energy distance. Namely, if $\XN^\circ$ is a critical point of the microscopic energy
\begin{equation}\label{lake:eq:ENdef}
\mathcal{H}_N(\ux_N) \coloneqq \frac{1}{2N}\sum_{1\leq i\neq j\leq N} \g(x_i-x_j) + \sum_{i=1}^N V(x_i),
\end{equation}
then  $\uz_N^t\coloneqq (x_i^t,v_i^t)_{i=1}^N$ with $(x_i^t, v_i^t) \coloneqq (x_i^\circ,0)$ is the unique (stationary) solution of \eqref{lake:eq:NewODE} with $\ga=0$, which is moreover \emph{independent of $\vep$}. In particular, if $\XN^\circ$ is a minimizer of $\mathcal{H}_N$, then by the results of \cite{sandier_2d_2015,sandier_1d_2015,rougerie_higher-dimensional_2016,petrache_next_2017} (see \cite[Chapters 11-12]{serfaty_lectures_2024} for a consolidated treatment), it holds that $\frac{1}{N}\sum_{i=1}^N \delta_{x_i^\circ} \xrightharpoonup[N\rightarrow\infty]{} \mu_V$ and
\begin{equation}
\frac{1}{N}\sum_{i=1}^N\underbrace{\zeta(x_i^\circ)}_{=0} + \Fr_N(\ux_N^\circ, \mu_V) + \frac{\log N}{2\d N}\indic_{\s=0} = \mathsf{C}_{\d,\s}^VN^{\frac{\s}{\d}-1} + o(N^{\frac{\s}{\d}-1}) \quad \text{as $N\rightarrow\infty$},
\end{equation}
where $\mathsf{C}_{\d,\s}^V$ is a computable constant depending only on $\d,\s,V$ {that encodes thermodynamic information at the microscale.} That $\zeta(x_i^\circ)=0$ for each $i$ follows from the fact that minimizing point configurations lie in $\supp\mu_V$, on which $\zeta$ vanishes.  Hence,
\begin{equation}
\Hr_{N}(\uz_N^t, 0)+  \frac{\log N}{2\d N\vep^2}\indic_{\s=0} = \frac{1}{\vep^2}\pa*{\Fr_N(\ux_N^\circ, \mu_V)+ \frac{\log N}{2\d N}\indic_{\s=0}} =  \mathsf{C}_{\d,\s}^V\frac{N^{\frac{\s}{\d}-1}}{\vep^2} + o(\frac{N^{\frac\s\d-1}}{\vep^2}),
\end{equation}
and we cannot expect vanishing of the total modulated energy in the supercritical mean-field regime if $\vep \leq N^{\frac{\s-\d}{2\d}}$.}

 In light of the preceding observations, one asks what is the effective equation describing the system \eqref{lake:eq:NewODE} as $\vep+\frac{1}{N}\rightarrow 0$, assuming that $\vep \leq N^{\frac{\s-\d}{2\d}}$. Naively, one might expect that when $\vep \ll N^{\frac{\s-\d}{2\d}}$, this behaves like first sending $\vep\rightarrow 0$ and then $N\rightarrow\infty$. However, this limit does not make sense in general. Indeed, multiplying both sides of the second equation of \eqref{lake:eq:NewODE} by $\vep^2$ and letting $\vep\rightarrow 0$ for fixed $N$, we formally see that the limiting positions  $\ux_{N}^t$ should be a critical point of the microscopic energy $\mathcal{H}_N$. If each limiting velocity $v_{i}^t$ is nonzero, we would not expect $\ux_{N}^t$ to remain a critical point for all $t$. On the other hand, there is no mechanism to force the velocities to tend to zero even if the initial positions are a critical point of the energy \eqref{lake:eq:ENdef}. In \cref{lake:ssec:1DCou}, we demonstrate this limit fails for the exactly solvable 1D Coulomb case precisely when $\vep^{-2}N^{\frac{\s}{\d}-1} = (\vep N)^{-2}$ does not vanish, showing that there need not be any weak limit for the empirical measure, even in the simplest setting.  


\begin{remark}\label{lake:rem:initMEvan}
Given $u^\circ$, one can produce statistically generic examples of initial data $\uz_N^\circ$ such that \eqref{lake:eq:SMFscl} holds. See \cite[Remark 1.2]{rosenzweig_lake_2025} for details.
\end{remark}

\subsection{Method of proof}\label{lake:ssec:introMPf}
We give some comments on the method of proof of \cref{lake:thm:mainSMF}. The complete details may be found in \cite[Section 5]{rosenzweig_lake_2025}.

{The quantity \eqref{lake:eq:soMEdef} is a variant of the total modulated energy originally introduced by Duerinckx and Serfaty \cite[Appendix]{serfaty_mean_2020} to treat the mean-field limit for Vlasov-Riesz in the monokinetic regime. The incorporation of the time-dependent corrector $\Uu^t$ in the modulated potential energy is inspired by \cite{han-kwan_newtons_2021}. The scaling by $\vep^2$ in  $\mu_V+\vep^2\Uu^t$ reflects $O(\vep^2)$ fluctuations around the macroscopic equilibrium spatial density $\mu_V$. The last term in \eqref{lake:eq:soMEdef} is a new contribution of \cite{rosenzweig_lake_2025} and reflects that the microscopic system \eqref{lake:eq:NewODE} is confined by an external potential $V$, as opposed to being set in a compact domain, e.g. $\T^\d$ as in \cite{han-kwan_newtons_2021, rosenzweig_rigorous_2023}. Although the $\zeta$ term appears to be only $O(1)$, it is in fact zero if the particles remain in the support of $\mu_V$, i.e. the quasineutral assumption is propagated. In analogy to the relationship between \cite{han-kwan_newtons_2021} and \cite{Brenier2000}, the total modulated energy \eqref{lake:eq:soMEdef} may be viewed as a renormalization of the total modulated energy from \cite{BCGM2015}, so as to allow for the Vlasov solution $f_\vep^t = \frac1N\sum_{i=1}^N \delta_{z_i^t}$.}

As with all modulated-energy approaches, the proof  is based on establishing a Gr\"{o}nwall relation for the total modulated energy \eqref{lake:eq:soMEdef}. The time derivative of this quantity has several terms that require different consideration. 

The main contribution from the modulated kinetic energy is trivially estimated using Cauchy-Schwarz. The correction $\vep^2\Uu^t$ in the spatial density is to cancel out the contribution of the pressure in \eqref{lake:eq:Lake} when one differentiates the modulated kinetic energy. 
The main contribution from the modulated potential energy is a first-order commutator \eqref{eq:introcomm} with $v=\tl{u}^t$, the extension of the solution $u^t$ of the Lake equation to the whole space, and $\mu=\mu_V+\vep^2\Uu^t$. To handle this term, we crucially rely on the sharp estimate \cref{thm:FIprime}\eqref{eq:FIsupC}, which allows for the scaling assumption \eqref{lake:eq:SMFscl}. 

A new term compared to \cite{han-kwan_newtons_2021, rosenzweig_rigorous_2023}, coming from the contribution of the $\zeta(x_i^t)$, is
\begin{align}\label{eq:zetaterm}
\frac1N\sum_{i=1}^N u^t(x_i^t)\cdot\nab\zeta(x_i^t) = \int_{\R^\d} u\cdot\nab\zeta d\mu_N^t,
\end{align}
This term has no commutator structure, but we instead manage to bound it by $\frac{C\|u^t\|_{W^{1,\infty}}}{N}\sum_{i=1}^N \zeta(x_i^t)$,
which is sufficient to close the Gr\"onwall loop. The proof of this bound uses that the no-flux condition  satisfied by $u^t$ on the boundary of $\Sigma = \supp\mu_V$ (a consequence of taking the quasineutral limit) and some nontrivial results for the regularity of the free boundary for solutions of the obstacle problem for the fractional Laplacian (see \cite[Appendix B]{rosenzweig_lake_2025}), which may be of independent interest. A term similar to \eqref{eq:zetaterm} was encountered in \cite{BCGM2015} in the Coulomb case, where the $\mu_N^t$ is replaced by the spatial density $\mu^t$ of Vlasov-Poisson, but handled by ad hoc arguments {using the local nature of the Laplacian}.

There are also several residual terms that have the form $\int_{\R^\d} \phi\,  d(\mu_N-\mu_V-\vep^2\Uu)$, where $\phi$ is a function of $u,\Uu$. These may be controlled by the modulated potential energy, thanks to its coercivity, plus errors which are $O(\vep^{-2}N^{\frac{\s}{\d}-1})$, hence acceptable.

Combining the estimates for the various terms and appealing to the Gr\"onwall-Bellman lemma yields the inequality \eqref{lake:eq:mainSMF}. The weak convergence of the empirical measure follows from the vanishing of the total modulated energy (see, e.g., \cite[Section 4.2]{rosenzweig_rigorous_2023}). This then completes the proof of \cref{lake:thm:mainSMF}.

\subsection{The regime $\vep^2 \lesssim N^{\frac{\s}{\d}-1}$ for the 1D Coulomb gas}\label{lake:ssec:1DCou}

Consider the 1D Coulomb case with quadratic confinement (i.e. $\s=-1$, $\g(x)= - |x|$, $V(x)=|x|^2$). Under the assumption that the initial empirical spatial density converges to the equilibrium measure, there is in general no weak limit for the empirical current as $\vep + N^{-1}\rightarrow 0$ if $\vep^{-2} N^{\frac{\s}{\d}-1} = \vep^{-2} N^{-2}$ does not vanish. Thus, for the supercritical mean-field limit in the (super-)Coulomb Riesz case, one cannot expect convergence to the Lake equation outside of the scaling assumption \eqref{lake:eq:SMFscl} in \cref{lake:thm:mainSMF}. This example uses the exact solvability of the one-dimensional Coulomb gas with quadratic confinement, which is folklore. We refer to \cite[Section 6]{rosenzweig_lake_2025} for the proof of the proposition stated below.

\begin{prop}\label{prop:1DCou}
There exists a sequence of initial positions $\XN^\circ$ with $\Fr_N(\XN^\circ,\mu_V)\rightarrow 0$  and $\mu_N^\circ \rightharpoonup \mu_V$ as $N\rightarrow\infty$, but that if {$v_i^\circ = 0$ for each $1\le i\le N$},  the associated empirical current $J_N^t \coloneqq \frac1N\sum_{i=1}^N v_i^t\delta_{x_i^t}$ 
\begin{itemize}
\item converges to zero (which is the unique solution of the Lake equation \eqref{lake:eq:LakeJ} {starting from zero initial datum}) uniformly on $[0,\infty)$ as $\vep+N^{-1}\rightarrow 0$  if $\vep N\rightarrow \infty$,
\item has no weak limit for any $t\in (0,\infty)$ as  $\vep+N^{-1}\rightarrow 0$ if $\vep N \not\rightarrow \infty$.
\end{itemize}
\end{prop}

\subsection{Sub-Coulomb case}\label{lake:ssec:subCou}
\cite{rosenzweig_lake_2025} only considered the (super-)Coulomb regime for two reasons. First, at the time of writing, only the nonsharp commutator estimate of \cite{nguyen_mean-field_2022} was available in the sub-Coulomb case, and so (super-)Coulomb was the only case where we could show our results were sharp. Obviously, this is no longer the case, thanks to \cref{thm:FIprime}\eqref{eq:FIsubC1}. 
Second, and more importantly, there does not seem to be an adequate regularity theory for the obstacle problem for higher powers of the fractional Laplacian (see \cite{DHP2023} for some progress in this direction). If one assumes that the equilibrium measure has full support in $\R^\d$, or restricts to the torus where full support is easily established under a smallness condition on the confinement, then all terms involving $\zeta$ vanish. One can then treat the sub-Coulomb case under the scaling assumption \eqref{lake:eq:SMFscl} by following the same proof in \cite{rosenzweig_lake_2025}.

{
\subsection{The non-monokinetic regime}\label{lake:ssec:nonmono}
The modulated-energy approach of \cite{rosenzweig_lake_2025} is currently limited to the monokinetic regime, as is the case for the kinetic mean-field limit. Given other approaches are capable of treating the mean-field limit for second-order systems without the monokinetic ansatz, a natural question is whether one can derive the full KIE equation \eqref{lake:eq:KIE} as $\vep + N^{-1}\rightarrow 0$, under some scaling relation between $\vep$ and $N$, drawing on these  approaches instead of the modulated-energy method.

The answer to this question is not immediately clear.  Per our knowledge, the only results deriving the full KIE (for the uniform density case) are due to Griffin-Pickering and Iacobelli \cite{GpI2018, GpI2020sing}, where one essentially sandwiches together two results: a mean-field limit for a (regularized) Vlasov-Poisson $N\rightarrow\infty$ and the quasineutral limit of Vlasov-Poisson. This approach imposes a very restrictive condition on the relation between $\vep$ and $N$.

At least in the case of regular interactions (treated in \cite[Appendix A]{rosenzweig_lake_2025}), preliminary calculations suggest that the modulated-energy approach of this paper can be extended beyond the monokinetic case, for possibly nonconstant density, through a multi-stream/multi-fluid decomposition, in the spirit of \cite{Grenier1996}. This direction is the subject of ongoing work. 
}

\bibliographystyle{alpha}
\bibliography{biblio}

\newcommand{\etalchar}[1]{$^{#1}$}
\begin{thebibliography}{CdCRS25}

\bibitem[AHS23]{altekruger_neural_2023}
Fabian Altekr{\"u}ger, Johannes Hertrich, and Gabriele Steidl.
\newblock Neural {Wasserstein} {Gradient} {Flows} for {Discrepancies} with {Riesz} {Kernels}.
\newblock In Andreas Krause, Emma Brunskill, Kyunghyun Cho, Barbara Engelhardt, Sivan Sabato, and Jonathan Scarlett, editors, {\em Proceedings of the 40th {International} {Conference} on {Machine} {Learning}}, volume 202 of {\em Proceedings of {Machine} {Learning} {Research}}, pages 664--690. PMLR, July 2023.

\bibitem[BBNY19]{BBNY2019}
Roland Bauerschmidt, Paul Bourgade, Miika Nikula, and Horng-Tzer Yau.
\newblock The two-dimensional {C}oulomb plasma: quasi-free approximation and central limit theorem.
\newblock {\em Adv. Theor. Math. Phys.}, 23(4):841--1002, 2019.

\bibitem[BCD11]{bahouri_fourier_2011}
Hajer Bahouri, Jean-Yves Chemin, and Rapha{\"e}l Danchin.
\newblock {\em Fourier analysis and nonlinear partial differential equations}, volume 343 of {\em Grundlehren der {Mathematischen} {Wissenschaften} [{Fundamental} {Principles} of {Mathematical} {Sciences}]}.
\newblock Springer, Heidelberg, 2011.

\bibitem[BCGM15]{BCGM2015}
Julien Barr\'{e}, David Chiron, Thierry Goudon, and Nader Masmoudi.
\newblock From {V}lasov-{P}oisson and {V}lasov-{P}oisson-{F}okker-{P}lanck systems to incompressible {E}uler equations: the case with finite charge.
\newblock {\em J. \'{E}c. polytech. Math.}, 2:247--296, 2015.

\bibitem[BDJ24]{BDJ2024}
Didier Bresch, Mitia Duerinckx, and Pierre-Emannuel Jabin.
\newblock A duality method for mean-field limits with singular interactions.
\newblock {\em arXiv preprint arXiv:2402.04695}, 2024.

\bibitem[BG80]{brezis_nonlinear_1980}
H.~Br\'ezis and T.~Gallouet.
\newblock Nonlinear {Schr{\"o}dinger} evolution equations.
\newblock {\em Nonlinear Anal.}, 4(4):677--681, 1980.
\newblock ISBN: 0362-546X, 1873-5215 Type: doi:10.1016/0362-546X(80)90068-1.

\bibitem[BG94]{BG1994}
Yann Brenier and Emmanuel Grenier.
\newblock Limite singuli\`ere du syst\`eme de {V}lasov-{P}oisson dans le r\'{e}gime de quasi neutralit\'{e}: le cas ind\'{e}pendant du temps.
\newblock {\em C. R. Acad. Sci. Paris S\'{e}r. I Math.}, 318(2):121--124, 1994.

\bibitem[BG13]{BG2013}
G.~Borot and A.~Guionnet.
\newblock Asymptotic expansion of {$\beta$} matrix models in the one-cut regime.
\newblock {\em Comm. Math. Phys.}, 317(2):447--483, 2013.

\bibitem[BG24]{BG2024}
Ga\"{e}tan Borot and Alice Guionnet.
\newblock Asymptotic expansion of {$\beta$} matrix models in the multi-cut regime.
\newblock {\em Forum Math. Sigma}, 12:Paper No. e13, 93, 2024.

\bibitem[BH77]{BH1977}
Walter Braun and Klaus Hepp.
\newblock The {V}lasov dynamics and its fluctuations in the {$1/N$} limit of interacting classical particles.
\newblock {\em Comm. Math. Phys.}, 56(2):101--113, 1977.

\bibitem[BHS19]{borodachov_discrete_2019}
Sergiy~V. Borodachov, Douglas~P. Hardin, and Edward~B. Saff.
\newblock {\em Discrete energy on rectifiable sets}.
\newblock Springer {Monographs} in {Mathematics}. Springer, New York, 2019.

\bibitem[BJS25]{BJS2022}
Didier Bresch, Pierre-Emmanuel Jabin, and Juan Soler.
\newblock A new approach to the mean-field limit of {V}lasov-{F}okker-{P}lanck equations.
\newblock {\em Anal. PDE}, 18(4):1037--1064, 2025.

\bibitem[BJW19a]{bresch_modulated_2019}
Didier Bresch, Pierre-Emmanuel Jabin, and Zhenfu Wang.
\newblock Modulated free energy and mean field limit.
\newblock {\em S\'eminaire Laurent Schwartza"EDP et applications}, pages 1--22, 2019.

\bibitem[BJW19b]{bresch_mean-field_2019}
Didier Bresch, Pierre-Emmanuel Jabin, and Zhenfu Wang.
\newblock On mean-field limits and quantitative estimates with a large class of singular kernels: application to the {Patlak}-{Keller}-{Segel} model.
\newblock {\em Comptes Rendus Math\'ematique. Acad\'emie des Sciences. Paris}, 357(9):708--720, 2019.

\bibitem[BJW23]{bresch_mean_2023}
Didier Bresch, Pierre-Emmanuel Jabin, and Zhenfu Wang.
\newblock Mean field limit and quantitative estimates with singular attractive kernels.
\newblock {\em Duke Mathematical Journal}, 172(13):2591--2641, 2023.

\bibitem[BLS18]{BLS2018}
Florent Bekerman, Thomas Lebl\'{e}, and Sylvia Serfaty.
\newblock C{LT} for fluctuations of {$\beta$}-ensembles with general potential.
\newblock {\em Electron. J. Probab.}, 23:Paper no. 115, 31, 2018.

\bibitem[B{\"O}19]{berman_propagation_2019}
Robert~J. Berman and Magnus {\"O}nnheim.
\newblock Propagation of {Chaos} for a {Class} of {First} {Order} {Models} with {Singular} {Mean} {Field} {Interactions}.
\newblock {\em SIAM Journal on Mathematical Analysis}, 51(1):159--196, January 2019.
\newblock Publisher: Society for Industrial and Applied Mathematics.

\bibitem[BP16]{BP2016}
Niklas Boers and Peter Pickl.
\newblock On mean field limits for dynamical systems.
\newblock {\em J. Stat. Phys.}, 164(1):1--16, 2016.

\bibitem[Bre89]{Brenier1989}
Yann Brenier.
\newblock {Une formulation de type Vlassov-Poisson pour les equations d'Euler des fluides parfaits incompressibles}.
\newblock Research Report RR-1070, 1989.

\bibitem[Bre00]{Brenier2000}
Y.~Brenier.
\newblock Convergence of the {V}lasov-{P}oisson system to the incompressible {E}uler equations.
\newblock {\em Comm. Partial Differential Equations}, 25(3-4):737--754, 2000.

\bibitem[BW80]{brezis_note_1980}
Haim Brezis and Stephen Wainger.
\newblock A note on limiting cases of sobolev embeddings and convolution inequalities.
\newblock {\em Communications in Partial Differential Equations}, 5(7):773--789, January 1980.
\newblock Publisher: Taylor \& Francis \_eprint: https://doi.org/10.1080/03605308008820154.

\bibitem[Cal80]{calderon_commutators_1980}
A.-P. Calder\'on.
\newblock Commutators, singular integrals on {Lipschitz} curves and applications.
\newblock In {\em Proceedings of the {International} {Congress} of {Mathematicians} ({Helsinki}, 1978)}, pages 85--96. Acad. Sci. Fennica, Helsinki, 1980.

\bibitem[CCH14]{carrillo_derivation_2014}
Jos\'e~Antonio Carrillo, Young-Pil Choi, and Maxime Hauray.
\newblock The derivation of swarming models: mean-field limit and {Wasserstein} distances.
\newblock In {\em Collective dynamics from bacteria to crowds}, volume 553 of {\em {CISM} {Courses} and {Lect}.}, pages 1--46. Springer, Vienna, 2014.

\bibitem[CD22]{CD2021}
Louis-Pierre Chaintron and Antoine Diez.
\newblock Propagation of chaos: a review of models, methods and applications. {I}. {M}odels and methods.
\newblock {\em Kinet. Relat. Models}, 15(6):895--1015, 2022.

\bibitem[CdCRS23]{chodron_de_courcel_sharp_2023}
Antonin Chodron~de Courcel, Matthew Rosenzweig, and Sylvia Serfaty.
\newblock Sharp uniform-in-time mean-field convergence for singular periodic {Riesz} flows.
\newblock {\em Annales de l'Institut Henri Poincar\'e C}, 42(2):391--472, December 2023.

\bibitem[CdCRS25]{chodron_de_courcel_attractive_2025}
Antonin Chodron~de Courcel, Matthew Rosenzweig, and Sylvia Serfaty.
\newblock The attractive log gas: {Stability}, uniqueness, and propagation of chaos.
\newblock {\em Communications of the American Mathematical Society}, 5:695--773, 2025.

\bibitem[CFGW24]{cai_propagation_2024}
Shuzhe Cai, Xuanrui Feng, Yun Gong, and Zhenfu Wang.
\newblock Propagation of chaos for 2d log gas on the whole space.
\newblock {\em arXiv preprint arXiv:2411.14777}, 2024.

\bibitem[CFP12]{carrillo_mass-transportation_2012}
Jos\'e~A. Carrillo, Lucas C.~F. Ferreira, and Juliana~C. Precioso.
\newblock A mass-transportation approach to a one dimensional fluid mechanics model with nonlocal velocity.
\newblock {\em Advances in Mathematics}, 231(1):306--327, 2012.

\bibitem[CJ87]{christ_polynomial_1987}
Michael Christ and Jean-Lin Journ\'e.
\newblock Polynomial growth estimates for multilinear singular integral operators.
\newblock {\em Acta Mathematica}, 159(1-2):51--80, 1987.

\bibitem[CM78]{coifman_commutateurs_1978}
R.~Coifman and Y.~Meyer.
\newblock Commutateurs d'int\'egrales singuli{\`e}res et op\'erateurs multilin\'eaires.
\newblock {\em Universit\'e de Grenoble. Annales de l'Institut Fourier}, 28(3):xi, 177--202, 1978.

\bibitem[CM22]{coraSpolyharmonicExtensionProblem2022}
Gabriele Cora and Roberta Musina.
\newblock The {\emph{s}}-polyharmonic extension problem and higher-order fractional {{Laplacians}}.
\newblock {\em Journal of Functional Analysis}, 283(5):109555, September 2022.

\bibitem[CN25]{cecchinConvergenceRateFluctuations2025}
Alekos Cecchin and Paul Nikolaev.
\newblock Convergence rate for {{Fluctuations}} of mean field interacting diffusion and application to {{2D}} viscous {{Vortex}} model and {{Coulomb}} potential.
\newblock {\em arXiv preprint arXiv:2509.01266}, 2025.

\bibitem[CR97]{CR1997}
S.~Jonathan Chapman and Giles Richardson.
\newblock Vortex pinning by inhomogeneities in type-ii superconductors.
\newblock {\em Physica D: Nonlinear Phenomena}, 108(4):397--407, 1997.

\bibitem[CS07]{CS2007}
Luis Caffarelli and Luis Silvestre.
\newblock An extension problem related to the fractional {L}aplacian.
\newblock {\em Comm. Partial Differential Equations}, 32(7-9):1245--1260, 2007.

\bibitem[DGR]{DGR}
Matias Delgadino, Rishabh Gvalani, and Matthew Rosenzweig.
\newblock Entropic commutator estimates.
\newblock Unpublished note.

\bibitem[DHAP23]{DHP2023}
Donatella Danielli, Alaa Haj~Ali, and Arshak Petrosyan.
\newblock The obstacle problem for a higher order fractional {L}aplacian.
\newblock {\em Calc. Var. Partial Differential Equations}, 62(8):Paper No. 218, 22, 2023.

\bibitem[DJ25]{duerinckxCorrelationEstimatesBrownian2025}
Mitia Duerinckx and Pierre-Emmanuel Jabin.
\newblock Correlation estimates for {{Brownian}} particles with singular interactions.
\newblock {\em arXiv preprint arXiv:2510.01507}, 2025.

\bibitem[Dob79]{Dobrushin1979}
R.~L. Dobru\v{s}in.
\newblock Vlasov equations.
\newblock {\em Funktsional. Anal. i Prilozhen.}, 13(2):48--58, 96, 1979.

\bibitem[DRAW02]{dauxois_dynamics_2002}
Thierry Dauxois, Stefano Ruffo, Ennio Arimondo, and Martin Wilkens.
\newblock Dynamics and thermodynamics of systems with long-range interactions: an introduction.
\newblock In {\em Dynamics and thermodynamics of systems with long-range interactions ({Les} {Houches}, 2002)}, volume 602 of {\em Lecture {Notes} in {Phys}.}, pages 1--19. Springer, Berlin, 2002.

\bibitem[DS18]{DS2018}
Mitia Duerinckx and Sylvia Serfaty.
\newblock Mean-field dynamics for {G}inzburg-{L}andau vortices with pinning and forcing.
\newblock {\em Ann. PDE}, 4(2):Paper No. 19, 172, 2018.

\bibitem[Due16]{duerinckx_mean-field_2016}
Mitia Duerinckx.
\newblock Mean-{Field} {Limits} for {Some} {Riesz} {Interaction} {Gradient} {Flows}.
\newblock {\em SIAM Journal on Mathematical Analysis}, 48(3):2269--2300, 2016.
\newblock \_eprint: https://doi.org/10.1137/15M1042620.

\bibitem[Due18]{Duerinckx2018}
Mitia Duerinckx.
\newblock Well-posedness for mean-field evolutions arising in superconductivity.
\newblock {\em Ann. Inst. H. Poincar\'{e} Anal. Non Lin\'{e}aire}, 35(5):1267--1319, 2018.
\newblock With an appendix jointly written with Julian Fischer.

\bibitem[Due21]{Duerinckx2021gl}
Mitia Duerinckx.
\newblock On the size of chaos via {G}lauber calculus in the classical mean-field dynamics.
\newblock {\em Comm. Math. Phys.}, 382(1):613--653, 2021.

\bibitem[Eng89]{engler_alternative_1989}
Hans Engler.
\newblock An {Alternative} {Proof} of the {Brezis} - {Wainger} {Inequality}.
\newblock {\em Communications in Partial Differential Equations}, 14(4):541--544, April 1989.
\newblock Publisher: Taylor \& Francis \_eprint: https://doi.org/10.1080/03605302.1989.12088448.

\bibitem[FKS82]{FKS1982}
Eugene~B. Fabes, Carlos~E. Kenig, and Raul~P. Serapioni.
\newblock The local regularity of solutions of degenerate elliptic equations.
\newblock {\em Comm. Partial Differential Equations}, 7(1):77--116, 1982.

\bibitem[FP25a]{feistl-heldMeanfieldLimitVlasovPoisson2025}
Manuela {Feistl-Held} and Peter Pickl.
\newblock On the mean-field limit for the {{Vlasov-Poisson}} system, April 2025.

\bibitem[FP25b]{feistl-heldMeanfieldLimitVlasovPoisson2025a}
Manuela {Feistl-Held} and Peter Pickl.
\newblock On the mean-field limit for the {{Vlasov-Poisson}} system in two dimensions, September 2025.

\bibitem[FS09]{FS2009}
Philippe Flajolet and Robert Sedgewick.
\newblock {\em Analytic combinatorics}.
\newblock Cambridge University Press, Cambridge, 2009.

\bibitem[FW24]{feng_quantitative_2024}
Xuanrui Feng and Zhenfu Wang.
\newblock Quantitative {Propagation} of {Chaos} for {2D} {Viscous} {Vortex} {Model} with {General} {Circulations} on the {Whole} {Space}.
\newblock {\em arXiv preprint arXiv:2411.14266}, 2024.

\bibitem[GBM24]{guillin_uniform_2024}
Arnaud Guillin, Pierre~Le Bris, and Pierre Monmarch\'e.
\newblock Uniform in time propagation of chaos for the {2D} vortex model and other singular stochastic systems.
\newblock {\em Journal of the European Mathematical Society}, 2024.

\bibitem[GBR{\etalchar{+}}06]{gretton_kernel_2006}
Arthur Gretton, Karsten Borgwardt, Malte Rasch, Bernhard Sch{\"o}lkopf, and Alex Smola.
\newblock A {Kernel} {Method} for the {Two}-{Sample}-{Problem}.
\newblock In B.~Sch{\"o}lkopf, J.~Platt, and T.~Hoffman, editors, {\em Advances in {Neural} {Information} {Processing} {Systems}}, volume~19. MIT Press, 2006.

\bibitem[GBR{\etalchar{+}}07]{gretton_kernel_2007}
Arthur Gretton, Karsten~M Borgwardt, Malte Rasch, Bernhard Sch{\"o}lkopf, and Alexander~J Smola.
\newblock A kernel approach to comparing distributions.
\newblock In {\em Proceedings of the national conference on artificial intelligence}, volume~22, page 1637. Menlo Park, CA; Cambridge, MA; London; AAAI Press; MIT Press; 1999, 2007.
\newblock Issue: 2.

\bibitem[GBR{\etalchar{+}}12]{gretton_kernel_2012}
Arthur Gretton, Karsten~M. Borgwardt, Malte~J. Rasch, Bernhard Sch{\"o}lkopf, and Alexander Smola.
\newblock A kernel two-sample test.
\newblock {\em Journal of Machine Learning Research (JMLR)}, 13:723--773, 2012.

\bibitem[Gol16]{golse_dynamics_2016}
Fran\c{c}ois Golse.
\newblock On the dynamics of large particle systems in the mean field limit.
\newblock In {\em Macroscopic and large scale phenomena: coarse graining, mean field limits and ergodicity}, volume~3 of {\em Lect. {Notes} {Appl}. {Math}. {Mech}.}, pages 1--144. Springer, 2016.

\bibitem[Gol22]{Golse2022ln}
Fran{\c{c}}ois Golse.
\newblock Mean-field limits in statistical dynamics.
\newblock {\em arXiv preprint arXiv:2201.02005}, 2022.

\bibitem[GP22]{GP2022}
Fran{\c c}ois Golse and Thierry Paul.
\newblock Mean-{{Field}} and {{Classical Limit}} for the {{N-Body Quantum Dynamics}} with {{Coulomb Interaction}}.
\newblock {\em Communications on Pure and Applied Mathematics}, 75(6):1332--1376, 2022.

\bibitem[GPI18]{GpI2018}
Megan Griffin-Pickering and Mikaela Iacobelli.
\newblock A mean field approach to the quasi-neutral limit for the {V}lasov-{P}oisson equation.
\newblock {\em SIAM J. Math. Anal.}, 50(5):5502--5536, 2018.

\bibitem[GPI20a]{GpI2020qn}
Megan Griffin-Pickering and Mikaela Iacobelli.
\newblock Recent developments on quasineutral limits for {V}lasov-type equations.
\newblock {\em arXiv preprint arXiv:2009.14169}, 2020.

\bibitem[GPI20b]{GpI2020sing}
Megan Griffin-Pickering and Mikaela Iacobelli.
\newblock Singular limits for plasmas with thermalised electrons.
\newblock {\em J. Math. Pures Appl. (9)}, 135:199--255, 2020.

\bibitem[Gra14]{grafakos_modern_2014}
Loukas Grafakos.
\newblock {\em Modern {Fourier} {Analysis}}.
\newblock Number 250 in Graduate {Texts} in {Mathematics}. Springer, third edition, 2014.

\bibitem[Gra21]{Grass2021}
Phillip Gra{\ss}.
\newblock Microscopic derivation of {Vlasov} equations with singular potentials.
\newblock {\em arXiv preprint arXiv:2105.06509}, 2021.

\bibitem[Gre95]{Grenier1995}
Emmanuel Grenier.
\newblock Defect measures of the {V}lasov-{P}oisson system in the quasineutral regime.
\newblock {\em Comm. Partial Differential Equations}, 20(7-8):1189--1215, 1995.

\bibitem[Gre96]{Grenier1996}
Emmanuel Grenier.
\newblock Oscillations in quasineutral plasmas.
\newblock {\em Comm. Partial Differential Equations}, 21(3-4):363--394, 1996.

\bibitem[Gre99]{Grenier1999}
Emmanuel Grenier.
\newblock Limite quasineutre en dimension 1.
\newblock In {\em Journ\'{e}es ``\'{E}quations aux {D}\'{e}riv\'{e}es {P}artielles'' ({S}aint-{J}ean-de-{M}onts, 1999)}, pages Exp. No. II, 8. Univ. Nantes, Nantes, 1999.

\bibitem[Han79]{hansson_imbedding_1979}
Kurt Hansson.
\newblock Imbedding theorems of {Sobolev} type in potential theory.
\newblock {\em Math. Scand.}, 45(1):77--102, 1979.
\newblock ISBN: 0025-5521, 1903-1807 Type: doi:10.7146/math.scand.a-11827.

\bibitem[Hau09]{hauray_wasserstein_2009}
Maxime Hauray.
\newblock Wasserstein distances for vortices approximation of {Euler}-type equations.
\newblock {\em Mathematical Models and Methods in Applied Sciences}, 19(8):1357--1384, 2009.

\bibitem[Hau14]{Hauray2014}
Maxime Hauray.
\newblock Mean field limit for the one dimensional {V}lasov-{P}oisson equation.
\newblock In {\em S\'{e}minaire {L}aurent {S}chwartz---\'{E}quations aux d\'{e}riv\'{e}es partielles et applications. {A}nn\'{e}e 2012--2013}, S\'{e}min. \'{E}qu. D\'{e}riv. Partielles, pages Exp. No. XXI, 16. \'{E}cole Polytech., Palaiseau, 2014.

\bibitem[HBGS23]{hertrich_wasserstein_2023}
Johannes Hertrich, Robert Beinert, Manuel Gr{\"a}f, and Gabriele Steidl.
\newblock Wasserstein gradient flows of the discrepancy with distance kernel on the line.
\newblock In Luca Calatroni, Marco Donatelli, Serena Morigi, Marco Prato, and Matteo Santacesaria, editors, {\em Scale {Space} and {Variational} {Methods} in {Computer} {Vision}}, pages 431--443, Cham, 2023. Springer International Publishing.

\bibitem[HC23]{hess-childs_large_2023}
Elias Hess-Childs.
\newblock Large deviation principles for singular {Riesz}-type diffusive flows.
\newblock {\em arXiv preprint arXiv:2312.02904}, 2023.

\bibitem[HCR25]{hCR2023}
Elias Hess-Childs and Keefer Rowan.
\newblock Higher-order propagation of chaos in {L}2 for interacting diffusions.
\newblock {\em Probab. Math. Phys.}, 6(2):581--646, 2025.

\bibitem[HCRS]{hess-childs_optimal_nodate}
Elias Hess-Childs, Matthew Rosenzweig, and Sylvia Serfaty.
\newblock Optimal quantization for {Riesz} {Maxmimum} {Mean} {Discrepancies}.
\newblock In preparation.

\bibitem[HCRS25]{hess-childs_sharp_2025}
Elias Hess-Childs, Matthew Rosenzweig, and Sylvia Serfaty.
\newblock A sharp commutator estimate for all {Riesz} modulated energies, November 2025.
\newblock arXiv:2511.13461 [math].

\bibitem[HHA{\etalchar{+}}23]{hagemann_posterior_2023}
Paul Hagemann, Johannes Hertrich, Fabian Altekr{\"u}ger, Robert Beinert, Jannis Chemseddine, and Gabriele Steidl.
\newblock Posterior sampling based on gradient flows of the {MMD} with negative distance kernel.
\newblock {\em arXiv preprint arXiv:2310.03054}, 2023.

\bibitem[HJ07]{HJ2007}
Maxime Hauray and Pierre-Emmanuel Jabin.
\newblock {$N$}-particles approximation of the {V}lasov equations with singular potential.
\newblock {\em Arch. Ration. Mech. Anal.}, 183(3):489--524, 2007.

\bibitem[HJ15]{HJ2015}
Maxime Hauray and Pierre-Emmanuel Jabin.
\newblock Particle approximation of {V}lasov equations with singular forces: propagation of chaos.
\newblock {\em Ann. Sci. \'{E}c. Norm. Sup\'{e}r. (4)}, 48(4):891--940, 2015.

\bibitem[HKH15]{HkH2015}
Daniel Han-Kwan and Maxime Hauray.
\newblock Stability issues in the quasineutral limit of the one-dimensional {V}lasov-{P}oisson equation.
\newblock {\em Comm. Math. Phys.}, 334(2):1101--1152, 2015.

\bibitem[HKI17a]{HkI2017jde}
Daniel Han-Kwan and Mikaela Iacobelli.
\newblock Quasineutral limit for {V}lasov-{P}oisson via {W}asserstein stability estimates in higher dimension.
\newblock {\em J. Differential Equations}, 263(1):1--25, 2017.

\bibitem[HKI17b]{HkI2017}
Daniel Han-Kwan and Mikaela Iacobelli.
\newblock The quasineutral limit of the {V}lasov-{P}oisson equation in {W}asserstein metric.
\newblock {\em Commun. Math. Sci.}, 15(2):481--509, 2017.

\bibitem[HKI21]{han-kwan_newtons_2021}
Daniel Han-Kwan and Mikaela Iacobelli.
\newblock From {Newton}'s second law to {Euler}'s equations of perfect fluids.
\newblock {\em Proceedings of the American Mathematical Society}, 149(7):3045--3061, 2021.

\bibitem[HKR16]{HkR2016}
Daniel Han-Kwan and Fr\'{e}d\'{e}ric Rousset.
\newblock Quasineutral limit for {V}lasov-{P}oisson with {P}enrose stable data.
\newblock {\em Ann. Sci. \'{E}c. Norm. Sup\'{e}r. (4)}, 49(6):1445--1495, 2016.

\bibitem[HLSS18]{hardin_large_2018}
Douglas~P. Hardin, Thomas Lebl\'e, Edward~B. Saff, and Sylvia Serfaty.
\newblock Large deviation principles for hypersingular {Riesz} gases.
\newblock {\em Constructive Approximation. An International Journal for Approximations and Expansions}, 48(1):61--100, 2018.

\bibitem[HM14]{hauray_kacs_2014}
Maxime Hauray and St\'ephane Mischler.
\newblock On {Kac}'s chaos and related problems.
\newblock {\em Journal of Functional Analysis}, 266(10):6055--6157, 2014.

\bibitem[HRS]{huang_fluctuations_nodate}
Jiaoyang Huang, Matthew Rosenzweig, and Sylvia Serfaty.
\newblock Fluctuations around the mean field limit for noisy singular {Riesz} flows.
\newblock In preparation.

\bibitem[HSSS17]{hardinNextOrderEnergy2017}
Douglas~P. Hardin, Edward~B. Saff, Brian~Z. Simanek, and Yujian Su.
\newblock Next order energy asymptotics for {{Riesz}} potentials on flat tori.
\newblock {\em International Mathematics Research Notices. IMRN}, (12):3529--3556, 2017.

\bibitem[HSST21]{hardin_dynamics_2021}
Douglas Hardin, Edward~B. Saff, Ruiwen Shu, and Eitan Tadmor.
\newblock Dynamics of particles on a curve with pairwise hyper-singular repulsion.
\newblock {\em Discrete and Continuous Dynamical Systems}, 41(12):5509--5536, 2021.

\bibitem[HWAH23]{hertrich_generative_2023}
Johannes Hertrich, Christian Wald, Fabian Altekr{\"u}ger, and Paul Hagemann.
\newblock Generative sliced {MMD} flows with {Riesz} kernels.
\newblock {\em arXiv preprint arXiv:2305.11463}, 2023.

\bibitem[Jab14]{Jabin2014}
Pierre-Emmanuel Jabin.
\newblock A review of the mean field limits for {V}lasov equations.
\newblock {\em Kinet. Relat. Models}, 7(4):661--711, 2014.

\bibitem[JW16]{JW2016}
Pierre-Emmanuel Jabin and Zhenfu Wang.
\newblock Mean field limit and propagation of chaos for {V}lasov systems with bounded forces.
\newblock {\em J. Funct. Anal.}, 271(12):3588--3627, 2016.

\bibitem[JW17]{JW2017_survey}
Pierre-Emmanuel Jabin and Zhenfu Wang.
\newblock {Mean field limit for stochastic particle systems}.
\newblock In {\em Act. Part. {V}ol. 1. {A}dvances theory, Model. Appl.}, Model. Simul. Sci. Eng. Technol., pages 379--402. Birkh{\"{a}}user/Springer, Cham, 2017.

\bibitem[JW18]{jabin_quantitative_2018}
Pierre-Emmanuel Jabin and Zhenfu Wang.
\newblock Quantitative estimates of propagation of chaos for stochastic systems with {W}$^{\textrm{-1,anfty}}$ kernels.
\newblock {\em Inventiones Mathematicae}, 214(1):523--591, 2018.

\bibitem[Kla10]{klainerman_pde_2010}
Sergiu Klainerman.
\newblock {PDE} as a {Unified} {Subject}.
\newblock In N.~Alon, J.~Bourgain, A.~Connes, M.~Gromov, and V.~Milman, editors, {\em Visions in {Mathematics}: {GAFA} 2000 {Special} {Volume}, {Part} {I}}, pages 279--315. Birkh{\"a}user, Basel, 2010.

\bibitem[KNSS22]{kolouri_generalized_2022}
Soheil Kolouri, Kimia Nadjahi, Shahin Shahrampour, and Umut Simsekli.
\newblock Generalized {Sliced} {Probability} {Metrics}.
\newblock In {\em {ICASSP} 2022 - 2022 {IEEE} {International} {Conference} on {Acoustics}, {Speech} and {Signal} {Processing} ({ICASSP})}, pages 4513--4517, 2022.

\bibitem[KP88]{KP1988}
Tosio Kato and Gustavo Ponce.
\newblock Commutator estimates and the {E}uler and {N}avier-{S}tokes equations.
\newblock {\em Comm. Pure Appl. Math.}, 41(7):891--907, 1988.

\bibitem[Lac23]{Lacker2023}
Daniel Lacker.
\newblock Hierarchies, entropy, and quantitative propagation of chaos for mean field diffusions.
\newblock {\em Probab. Math. Phys.}, 4(2):377--432, 2023.

\bibitem[Laz16]{Lazarovici2016}
Dustin Lazarovici.
\newblock The {V}lasov-{P}oisson dynamics as the mean field limit of extended charges.
\newblock {\em Comm. Math. Phys.}, 347(1):271--289, 2016.

\bibitem[Li19]{Li2019}
Dong Li.
\newblock On {K}ato-{P}once and fractional {L}eibniz.
\newblock {\em Rev. Mat. Iberoam.}, 35(1):23--100, 2019.

\bibitem[LLF23]{LlF2023}
Daniel Lacker and Luc Le~Flem.
\newblock Sharp uniform-in-time propagation of chaos.
\newblock {\em Probability Theory and Related Fields}, 2023.

\bibitem[LLN20]{LLN2020}
Tau~Shean Lim, Yulong Lu, and James~H. Nolen.
\newblock Quantitative propagation of chaos in a bimolecular chemical reaction-diffusion model.
\newblock {\em SIAM J. Math. Anal.}, 52(2):2098--2133, 2020.

\bibitem[LOT96a]{LOT1996}
C.~David Levermore, Marcel Oliver, and Edriss~S. Titi.
\newblock Global well-posedness for models of shallow water in a basin with a varying bottom.
\newblock {\em Indiana Univ. Math. J.}, 45(2):479--510, 1996.

\bibitem[LOT96b]{LOT1996phy}
C.David Levermore, Marcel Oliver, and Edriss~S. Titi.
\newblock Global well-posedness for the lake equations.
\newblock {\em Physica D: Nonlinear Phenomena}, 98(2):492--509, 1996.
\newblock Nonlinear Phenomena in Ocean Dynamics.

\bibitem[LP17]{LP2017}
Dustin Lazarovici and Peter Pickl.
\newblock A mean field limit for the {V}lasov-{P}oisson system.
\newblock {\em Arch. Ration. Mech. Anal.}, 225(3):1201--1231, 2017.

\bibitem[LS18]{leble_fluctuations_2018}
Thomas Lebl\'e and Sylvia Serfaty.
\newblock Fluctuations of two dimensional {Coulomb} gases.
\newblock {\em Geometric and Functional Analysis}, 28(2):443--508, 2018.

\bibitem[M\'24]{menard_mean-field_2024}
Matthieu M\'enard.
\newblock Mean-field limit derivation of a monokinetic spray model with gyroscopic effects.
\newblock {\em SIAM Journal on Mathematical Analysis}, 56(1):1068--1113, 2024.

\bibitem[Mas01]{Masmoudi2001}
Nader Masmoudi.
\newblock From {V}lasov-{P}oisson system to the incompressible {E}uler system.
\newblock {\em Comm. Partial Differential Equations}, 26(9-10):1913--1928, 2001.

\bibitem[Mas07]{Masmoudi2007}
Nader Masmoudi.
\newblock Rigorous derivation of the anelastic approximation.
\newblock {\em J. Math. Pures Appl. (9)}, 88(3):230--240, 2007.

\bibitem[MD24]{modeste_characterization_2024}
Thibault Modeste and Cl\'ement Dombry.
\newblock Characterization of translation invariant {MMD} on \{\${\textbackslash}{bbb} {r}{\textasciicircum}d\$\} and connections with {Wasserstein} distances.
\newblock {\em Journal of Machine Learning Research (JMLR)}, 25:Paper No. [237], 39, 2024.

\bibitem[M{\'e}n23]{Menard2023lake}
Matthieu M{\'e}nard.
\newblock Mean-field limit of point vortices for the lake equations.
\newblock {\em arXiv preprint arXiv:2309.10453}, 2023.

\bibitem[MKTL08]{MKTL2008}
JT~Mendon{\c{c}}a, R~Kaiser, H~Ter{\c{c}}as, and J~Loureiro.
\newblock Collective oscillations in ultracold atomic gas.
\newblock {\em Phys. Rev. A}, 78:013408, Jul 2008.

\bibitem[Mos70]{moser_sharp_1970}
J.~Moser.
\newblock A sharp form of an inequality by {N}. {Trudinger}.
\newblock {\em Indiana University Mathematics Journal}, 20:1077--1092, 1970.

\bibitem[M{\"u}l97]{muller_integral_1997}
Alfred M{\"u}ller.
\newblock Integral probability metrics and their generating classes of functions.
\newblock {\em Advances in Applied Probability}, 29(2):429--443, 1997.

\bibitem[NRS22]{nguyen_mean-field_2022}
Quoc-Hung Nguyen, Matthew Rosenzweig, and Sylvia Serfaty.
\newblock Mean-field limits of {Riesz}-type singular flows.
\newblock {\em Ars Inveniendi Analytica}, pages Paper No. 4, 45, 2022.

\bibitem[NW74]{NW1974}
H.~Neunzert and J.~Wick.
\newblock Die {A}pproximation der {L}\"{o}sung von {I}ntegro-{D}ifferentialgleichungen durch endliche {P}unktmengen.
\newblock In {\em Numerische {B}ehandlung nichtlinearer {I}ntegrodifferential- und {D}ifferentialgleichungen ({T}agung, {M}ath. {F}orschungsinst., {O}berwolfach, 1973)}, pages 275--290. Lecture Notes in Math., Vol. 395. 1974.

\bibitem[OD98]{OD1998}
T.~M. O'Neil and Daniel H.~E. Dubin.
\newblock Thermal equilibria and thermodynamics of trapped plasmas with a single sign of charge.
\newblock {\em Physics of Plasmas}, 5(6):2163--2193, 1998.

\bibitem[OP62]{OP1962}
Yoshimitsu Ogura and Norman~A. Phillips.
\newblock Scale analysis of deep and shallow convection in the atmosphere.
\newblock {\em Journal of Atmospheric Sciences}, 19(2):173 -- 179, 1962.

\bibitem[PCJ25]{porat_singular_2025}
Immanuel~Ben Porat, Jos\'e~A Carrillo, and Pierre-Emmanuel Jabin.
\newblock Singular flows with time-varying weights.
\newblock {\em arXiv preprint arXiv:2503.02276}, 2025.

\bibitem[Por23]{porat_derivation_2023}
Immanuel~Ben Porat.
\newblock Derivation of {Euler}'s equations of perfect fluids from von {Neumann}'s equation with magnetic field.
\newblock {\em Journal of Statistical Physics}, 190(7):Paper No. 121, 44, 2023.

\bibitem[PS17]{petrache_next_2017}
Mircea Petrache and Sylvia Serfaty.
\newblock Next order asymptotics and renormalized energy for {Riesz} interactions.
\newblock {\em Journal of the Institute of Mathematics of Jussieu. JIMJ. Journal de l'Institut de Math\'ematiques de Jussieu}, 16(3):501--569, 2017.

\bibitem[PS25]{peilen_local_2025}
Luke Peilen and Sylvia Serfaty.
\newblock Local {Laws} and {Fluctuations} for {Super}-{Coulombic} {Riesz} {Gases}, December 2025.
\newblock arXiv:2511.18623 [math-ph].

\bibitem[Ric00]{Richardson2000}
G.~Richardson.
\newblock Vortex motion in shallow water with varying bottom topography and zero {F}roude number.
\newblock {\em J. Fluid Mech.}, 411:351--374, 2000.

\bibitem[Ros]{rosenzweig_cumulants_nodate}
Matthew Rosenzweig.
\newblock Cumulants of {Riesz} gases.
\newblock Unpublished note.

\bibitem[Ros20]{rosenzweig_mean-field_2020}
Matthew Rosenzweig.
\newblock The {Mean}-{Field} {Limit} of {Stochastic} {Point} {Vortex} {Systems} with {Multiplicative} {Noise}.
\newblock {\em arXiv preprint arXiv:2011.12180}, 2020.
\newblock Forthcoming in Comm. Pure Appl. Math.

\bibitem[Ros21]{rosenzweig_quantum_2021}
Matthew Rosenzweig.
\newblock From {Quantum} {Many}-{Body} {Systems} to {Ideal} {Fluids}.
\newblock {\em arXiv preprint arXiv:2110.04195}, 2021.

\bibitem[Ros22a]{rosenzweigMeanfieldApproximationHigherdimensional2022}
Matthew Rosenzweig.
\newblock The mean-field approximation for higher-dimensional {{Coulomb}} flows in the scaling-critical {{L}}\textasciicircum\i nfty space.
\newblock {\em Nonlinearity}, 35(6):2722--2766, May 2022.

\bibitem[Ros22b]{rosenzweig_mean-field_2022-1}
Matthew Rosenzweig.
\newblock Mean-{Field} {Convergence} of {Point} {Vortices} to the {Incompressible} {Euler} {Equation} with {Vorticity} in {L}{\textasciicircum}anfty.
\newblock {\em Archive for Rational Mechanics and Analysis}, 243(3):1361--1431, 2022.

\bibitem[Ros23]{rosenzweig_rigorous_2023}
Matthew Rosenzweig.
\newblock On the rigorous derivation of the incompressible {Euler} equation from {Newton}'s second law.
\newblock {\em Letters in Mathematical Physics}, 113(1):Paper No. 13, 32, 2023.

\bibitem[RSa]{rosenzweig_commutator_nodate}
Matthew Rosenzweig and Sylvia Serfaty.
\newblock Commutator estimates, {Stein}'s method, and the transport approach to fluctuations of {Riesz} gases.
\newblock In preparation.

\bibitem[RSb]{rosenzweig_sharp_nodate}
Matthew Rosenzweig and Sylvia Serfaty.
\newblock Sharp commutator estimates of all order for {Coulomb} and {Riesz} modulated energies.
\newblock {\em Communications on Pure and Applied Mathematics}, n/a(n/a).
\newblock \_eprint: https://onlinelibrary.wiley.com/doi/pdf/10.1002/cpa.70010.

\bibitem[RS16]{rougerie_higher-dimensional_2016}
Nicolas Rougerie and Sylvia Serfaty.
\newblock Higher-dimensional {Coulomb} gases and renormalized energy functionals.
\newblock {\em Communications on Pure and Applied Mathematics}, 69(3):519--605, 2016.

\bibitem[RS23]{rosenzweig_global--time_2023}
Matthew Rosenzweig and Sylvia Serfaty.
\newblock Global-in-time mean-field convergence for singular {Riesz}-type diffusive flows.
\newblock {\em The Annals of Applied Probability}, 33(2):754--798, 2023.

\bibitem[RS24]{rosenzweig_relative_2024}
Matthew Rosenzweig and Sylvia Serfaty.
\newblock Relative entropy and modulated free energy without confinement via self-similar transformation.
\newblock {\em arXiv preprint arXiv:2402.13977}, 2024.

\bibitem[RS25a]{rosenzweig_lake_2025}
Matthew Rosenzweig and Sylvia Serfaty.
\newblock The lake equation as a supercritical mean-field limit.
\newblock {\em Journal de l'\'Ecole polytechnique. Math\'ematiques}, 12:1019--1068, 2025.

\bibitem[RS25b]{rosenzweig_modulated_2025}
Matthew Rosenzweig and Sylvia Serfaty.
\newblock Modulated logarithmic {Sobolev} inequalities and generation of chaos.
\newblock {\em Annales de la Facult\'e des Sciences de Toulouse. Math\'ematiques. S\'erie 6}, 34(1):107--134, 2025.

\bibitem[RSW]{rosenzweig_wasserstein_nodate}
Matthew Rosenzweig, Dejan Slep{\v{c}}ev, and Lihan Wang.
\newblock Wasserstein gradient flow of {Maximum} {Mean} {Discrepancy} with energy kernels.

\bibitem[Rub24]{Rubin2024}
Boris Rubin.
\newblock {\em Fractional integrals, potentials, and {R}adon transforms}.
\newblock Chapman and Hall/CRC, 2024.

\bibitem[Ser20]{serfaty_mean_2020}
Sylvia Serfaty.
\newblock Mean field limit for {Coulomb}-type flows.
\newblock {\em Duke Mathematical Journal}, 169(15):2887--2935, October 2020.
\newblock Publisher: Duke University Press.

\bibitem[Ser23]{serfaty_gaussian_2023}
Sylvia Serfaty.
\newblock Gaussian fluctuations and free energy expansion for {Coulomb} gases at any temperature.
\newblock {\em Annales de l'Institut Henri Poincar\'e Probabilit\'es et Statistiques}, 59(2):1074--1142, 2023.

\bibitem[Ser24]{serfaty_lectures_2024}
Sylvia Serfaty.
\newblock Lectures on {Coulomb} and {Riesz} gases.
\newblock {\em arXiv preprint arXiv:2407.21194}, 2024.

\bibitem[SS15a]{sandier_1d_2015}
Etienne Sandier and Sylvia Serfaty.
\newblock {1D} log gases and the renormalized energy: crystallization at vanishing temperature.
\newblock {\em Probability Theory and Related Fields}, 162(3-4):795--846, 2015.

\bibitem[SS15b]{sandier_2d_2015}
Etienne Sandier and Sylvia Serfaty.
\newblock {2D} {Coulomb} gases and the renormalized energy.
\newblock {\em The Annals of Probability}, 43(4):2026--2083, 2015.

\bibitem[SSS19]{seeger_multilinear_2019}
Andreas Seeger, Charles~K. Smart, and Brian Street.
\newblock Multilinear singular integral forms of {Christ}-{Journ\'e} type.
\newblock {\em Memoirs of the American Mathematical Society}, 257(1231):v+134, 2019.
\newblock ISBN: 978-1-4704-3437-3; 978-1-4704-4945-2.

\bibitem[Ste70]{stein_singular_1970}
Elias~M Stein.
\newblock {\em Singular {Integrals} and {Differentiability} {Properties} of {Functions}}, volume~2.
\newblock Princeton University Press, 1970.

\bibitem[Szn91]{Sznitman1991}
Alain-Sol Sznitman.
\newblock Topics in propagation of chaos.
\newblock In {\em \'{E}cole d'\'{E}t\'{e} de {P}robabilit\'{e}s de {S}aint-{F}lour {XIX}---1989}, volume 1464 of {\em Lecture Notes in Math.}, pages 165--251. Springer, Berlin, 1991.

\bibitem[Tro86]{Trocheris1986}
M.~Trocheris.
\newblock On the derivation of the one-dimensional {V}lasov equation.
\newblock {\em Transport Theory Statist. Phys.}, 15(5):597--628, 1986.

\bibitem[Tru67]{trudinger_imbeddings_1967}
Neil~S. Trudinger.
\newblock On imbeddings into {Orlicz} spaces and some applications.
\newblock {\em J. Math. Mech.}, 17:473--483, 1967.

\bibitem[Uhl18]{Uhlemann2018}
Cora Uhlemann.
\newblock Finding closure: approximating {Vlasov-Poisson} using finitely generated cumulants.
\newblock 2018(10):030, oct 2018.

\bibitem[Wan26]{Wang2024sharp}
Songbo Wang.
\newblock Sharp local propagation of chaos for mean field particles with {$W^{-1,\infty}$} kernels.
\newblock {\em J. Funct. Anal.}, 290(3):Paper No. 111240, 2026.

\bibitem[WBIP85]{WBIP1985}
D.~J. Wineland, J.~J. Bollinger, Wayne~M. Itano, and J.~D. Prestage.
\newblock Angular momentum of trapped atomic particles.
\newblock {\em J. Opt. Soc. Am. B}, 2(11):1721--1730, Nov 1985.

\bibitem[Yan13]{yangHigherOrderExtensions2013}
Ray Yang.
\newblock On higher order extensions for the fractional {{Laplacian}}, February 2013.

\end{thebibliography}
\end{document}